\documentclass[journal,article,submit,moreauthors,12pt,a4paper]{mdpi} 
\lastpage{x}
\doinum{10.3390/------}
\pubvolume{xx}
\pubyear{2013}
\history{Received: xx / Accepted: xx / Published: xx}

\usepackage{amsmath, amssymb, amsthm}
\usepackage{latexsym}
\usepackage{indentfirst}
\usepackage{graphicx}
\usepackage{placeins}
\usepackage{booktabs}
\usepackage{algorithm}
\usepackage{algorithmic}
\usepackage{subfigure}

\theoremstyle{plain}

\theoremstyle{definition}

\theoremstyle{remark}
\newtheorem*{remark}{Remark}

\newcommand{\videpost}{\textit{vide post}{}}
\newcommand{\etc}{\textit{etc}{}}

\newcommand{\ie}{\textit{i.e.}{~}}
\newcommand{\eg}{\textit{e.g.}{~}}

\newcommand{\rd}{\mathrm{d}}
\newcommand{\ud}{\,\mathrm{d}}

\newcommand{\Or}{\mathcal{O}}

\newcommand{\bd}[1]{\boldsymbol{#1}}
\newcommand{\wt}[1]{\widetilde{#1}}

\newcommand{\wb}[1]{\overline{#1}}

\newcommand{\mc}[1]{\mathcal{#1}}

\newcommand{\mcF}[1]{\mathcal{F}}

\newcommand{\abs}[1]{\left\lvert#1\right\rvert}

\newcommand{\ion}{\mathrm{ion}}
\newcommand{\hxc}{\mathrm{hxc}}

\newcommand{\rhoscf}{\rho_{\mathrm{SCF}}}
\newcommand{\err}{\mathrm{err}}
\newcommand{\E}{\mathrm{E}}

\newcommand{\barint}{\kern4pt \raise3.4pt\hbox{\vrule height.6pt
    width7pt} \kern-11pt \int}

\newcommand{\vL}{\mathbf{\mathcal{L}}}
\newcommand{\vR}{\mathbf{R}}
\newcommand{\vv}{\mathbf{v}}
\newcommand{\vzero}{\mathbf{0}}
\newcommand{\Ker}{\mathrm{Ker}}

\Title{Analysis of the Time Reversible Born-Oppenheimer Molecular Dynamics}

\Author{Lin Lin $^{1}$, Jianfeng Lu $^{2}$ and Sihong Shao $^{3,\ast}$}

\address{%
$^{1}$ Computational Research Division, Lawrence Berkeley National
Laboratory, Berkeley, CA 94720, USA. Email: \url{linlin@lbl.gov}\\
$^{2}$ Department of Mathematics and Department of Physics, Duke University, Box 90320, Durham, NC 27708, USA. Email:  \url{jianfeng@math.duke.edu} \\
$^{3}$ LMAM and School of Mathematical Sciences, Peking University,
Beijing 100871, China. Email: \url{sihong@math.pku.edu.cn}}

\corres{Tel: +86-10-62753433, Fax: +86-10-62751801.}

\abstract{We analyze the time reversible Born-Oppenheimer molecular
dynamics (TRBOMD) scheme, which preserves the time reversibility of the
Born-Oppenheimer molecular dynamics even with non-convergent
self-consistent field iteration.  In the linear response regime,
we derive the stability condition as well as the accuracy of TRBOMD for
computing physical properties such as the phonon frequency obtained from
the molecular dynamic simulation. We connect and compare TRBOMD with
the Car-Parrinello molecular dynamics in terms of accuracy and
stability. We further discuss the accuracy of TRBOMD
beyond the linear response regime for non-equilibrium dynamics of nuclei.
Our results are demonstrated through numerical experiments using a
simplified one dimensional model for Kohn-Sham density functional
theory.}



\begin{document}

\section{Introduction}\label{sec:intro}

\textit{Ab initio} molecular dynamics
(AIMD)~\cite{MarxHutter2000,Kirchner2012,PayneTeterAllenEtAl1992,DeumensDizLongoEtAl1994,TuckermanUngarRosenvingeEtAl1996,Parrinello1997}
has been greatly developed
in the past few decades, so that nowadays it is able to quantitatively
predict the equilibrium and non-equilibrium properties {for} a vast
range of systems. {AIMD has} become widely used in chemistry, biology,
materials science \etc.  Most AIMD methods treat the nuclei as
classical particles following the Newtonian dynamics ({known as the
  time dependent} Born-Oppenheimer approximation), and the interactive
force among nuclei is provided directly from electronic structure
theory, such as the Kohn-Sham density functional
theory~\cite{HohenbergKohn1964,KohnSham1965} (KSDFT), without the need
of {using empirical atomic potentials}.  KSDFT consists of a set of
nonlinear equations {which are solved at each molecular dynamics time
  step} \textit{self-consistently} via the self-consistent field (SCF)
iteration.  In the Born-Oppenheimer molecular dynamics (BOMD), KSDFT
is solved till full self-consistency for each atomic configuration per
time step.  Since many iterations are usually needed to reach full
self-consistency and each iteration takes considerable amount of time,
{until recently} this procedure was still found to be prohibitively
expensive for producing meaningful dynamical information. On the other
hand, if the self-consistent iterations are truncated {before
  convergence} is reached, it is often the case that the energy of the
system {is no longer conservative} even for an NVE system. The error
in SCF iteration acts as a sink or source, gradually draining or
adding energy to the atomic system within a short period of molecular
dynamics simulation~\cite{RemlerMadden1990}. This is one of the main challenges for accelerating Born-Oppenheimer molecular dynamics.

AIMD was made practical by the
ground-breaking work of Car-Parrinello molecular dynamics
(CPMD)~\cite{CarParrinello1985}.  CPMD introduces an extended
Lagrangian including the degrees of freedom of both nuclei and
electrons without the necessity of a convergent SCF iteration.  The
dynamics of electronic orbitals can be loosely viewed as a special way
for performing the SCF iteration at each molecular dynamics (MD) step.
Thanks to the Hamiltonian structure, numerical simulation for CPMD is
stable, and the energy is conservative over a much longer time {period} compared to
that for BOMD with non-convergent SCF iteration.  {When the system
has a spectral gap, the accuracy of CPMD is controlled by a single
parameter, the fictitious electron mass $\mu$. The result of CPMD
approaches that of BOMD as $\mu$ goes to
zero~\cite{PastoreSmargiassiBuda1991, BornemannSchutte1998}.  However,
it has also been shown that CPMD does not work as well for systems
with vanishing gap, for example for metallic
systems~\cite{PastoreSmargiassiBuda1991}.

To reduce the cost of BOMD, in particular, the number of SCF
iterations needed per MD time step, a new type of AIMD method, the
time reversible Born-Oppenheimer molecular dynamics (TRBOMD) method
has been recently proposed by Niklasson, Tymczak and Challacombe
in~\cite{Niklasson2006}. The method has been further developed
in~\cite{Niklasson2007,Niklasson2008,Niklasson2009,Niklasson2012}.
The idea of TRBOMD can be summarized as follows: TRBOMD assumes that
the SCF iteration is a \textit{deterministic} procedure,
with the outcome determined only by the initial guess of the variable
to be determined self-consistently.  For instance, this variable can
be the electron density, and the SCF iteration procedure
can be simple mixing with a fixed number of iteration steps without
reaching full self-consistency. Then a fictitious dynamics governed by
a second order ordinary differential equation (ODE) is introduced on
this initial guess variable.  The resulting coupled dynamics is then
time-reversible and supposed to be more stable since it has been found
that time-reversible numerical schemes are more stable for long time
simulation~\cite{HairerLubichWanner, McLachlan2004}.

Although TRBOMD has been found to be effective and significantly reduces
the number of {SCF} iterations needed in practice, to the
extent of our knowledge there has been so far no detailed analysis
of TRBOMD, other than the numerical stability condition of
the Verlet or generalized Verlet scheme for time
discretization~\cite{Niklasson2009}. Accuracy, stability, as well as
the applicability range of TRBOMD remain unclear. In particular, it is
not known how the choice of SCF iteration scheme affects TRBOMD}. These are crucial issues for guiding the practical {use} of
TRBOMD.  The full TRBOMD
method for general systems is highly nonlinear and is difficult to
analyze.  In this work, we first focus on the linear response regime,
\ie we assume that each atom oscillates around their equilibrium
position and the electron density stays around the ``true'' electron
density. Under such assumptions, we analyze the
accuracy and stability of TRBOMD. We then extend the results to the regime where the atom position is not near equilibrium using averaging principle.

The rest of the paper is organized as follows. We illustrate the idea of
TRBOMD and its analysis in the linear response regime using a
simple model in Section~\ref{sec:observation}, and introduce TRBOMD
for AIMD in Section~\ref{sec:trbomd}. We analyze TRBOMD in
the linear response regime, and compare TRBOMD with CPMD in
Section~\ref{sec:lr}.  The numerical results for TRBOMD in the
linear response regime are given in Section~\ref{sec:numer}.  We present
the analysis of TRBOMD beyond the linear response regime such as the
non-equilibrium dynamics in Section~\ref{sec:nonlinear},
and conclude with a few remarks in
Section~\ref{sec:conclusion}.

\section{An illustrative model}\label{sec:observation}

To start, let us illustrate the main idea for a simple model problem, {which provides the essence of TRBOMD in a much simplified setting}.
Consider the following nonlinear ODE
\begin{equation}
	\ddot{x}(t) = f(x(t))
	\label{eqn:ODE}
\end{equation}
{where we assume} that the {right hand side} $f(x)$ is difficult to compute, {and it can be approximated by an iterative procedure}.  Starting from {an} initial guess
$s\approx f(x)$, the {final approximation via the iterative procedure is denoted by $g(x,s)$}.  {We assume the
approximation} $g(x,s)$ {is consistent}, \ie
\begin{equation}
	g(x,f(x)) = f(x).
	\label{eqn:gconsistency}
\end{equation}
{To numerically solve the ODE~\eqref{eqn:ODE}, we discretize it by some numerical scheme, then it remains to decide the initial guess $s$ at each time step.} {A natural
choice of $s$ would be $g(x,s)$ from the previous step, as $x$ does not change much in successive steps}.  For instance, if the
Verlet algorithm is used and $t_{k}=k\Delta t$ with $\Delta t$ being the time step, the discretized ODE
becomes
\begin{equation}
	\begin{split}
	x_{k+1} &= 2 x_{k} - x_{k-1} + (\Delta t)^2 g(x_{k}, s_{k}), \\
	s_{k+1} &= g(x_{k}, s_{k}).
	\end{split}
	\label{eqn:VerletODE}
\end{equation}
{We immediately observe that the discretization
scheme~\eqref{eqn:VerletODE} breaks the time reversibility of the
original ODE~\eqref{eqn:ODE}. In other words, for the original
ODE~\eqref{eqn:ODE}, we propagate the system forward in time from
$(x(t_0), \dot{x}(t_0))$ to $(x(t_1), \dot{x}(t_1))$. Then if we use
$(x(t_1), \dot{x}(t_1))$ as the initial data at $t = t_1$ and propagate
the system backward in time to time $t = t_0$, we will be at the state
$(x(t_0), \dot{x}(t_0))$. The loss of the time reversible structure can
introduce large error in long time numerical
simulation~\cite{McLachlan2004}.  This is the main reason why BOMD with
non-convergent SCF iteration fails for long time
simulations~\cite{Niklasson2006}.}
To overcome this obstacle, {the idea of TRBOMD is to
introduce a fictitious dynamics for the initial guess $s$.
Namely, we consider the time reversible coupled system}
\begin{equation}
	\begin{split}
		\ddot{x}(t) &= g(x(t), s(t)),\\
		\ddot{s}(t) &= \omega^2 (g(x(t),s(t))-s(t)),
	\end{split}
	\label{eqn:TRBOMD_ODE}
\end{equation}
where $\omega$ is an artificial frequency.
We analyze now the accuracy and stability of Eq.~\eqref{eqn:TRBOMD_ODE}
in the linear response regime by assuming that the trajectory $x(t)$
oscillates around a equilibrium position $x^{*}$.
We denote by $\wt{x}(t)=x(t)-x^{*}$ the deviation from the equilibrium position
and $\wt{s}(t)=s(t)-f(x(t))$ the deviation of the initial guess
from the exact force term.
Consequently, the equation of
motion~\eqref{eqn:TRBOMD_ODE} can be rewritten as (for simplicity we
suppress the $t$-dependence in the notation for the rest of the section)
\begin{equation}
	\begin{split}
		\ddot{\wt{x}} & = g(x,s),\\
		\ddot{\wt{s}} & = \omega^2 (g(x,s) - s) - f''(x)
		(\dot{x})^2 - f'(x) \ddot{x}.
	\end{split}
	\label{eqn:TRBOMD_ODE1}
\end{equation}
where the term $- f''(x)
		(\dot{x})^2 - f'(x) \ddot{x}$ comes from the
term $f(x)$ in $\wt{s}$ by the chain rule.

In the linear response regime, we assume {the linear approximation of force for $x$ around $x^*$}:
\begin{equation}
	f(x) \approx -\Omega^2 (x-x^*) = -\Omega^2 \wt{x},
	\label{eqn:flinear}
\end{equation}
where $\Omega$ is the oscillation frequency of $x$ in the linear response
regime.  We {also} linearize {$g$} with respect to $\wt{s}$ and $\wt{x}$
and dropping all higher order terms as
\begin{equation}
	\begin{split}
		g(x,s) &= g(x,f(x)+\wt{s}) \\
		&\approx g(x,f(x)) + g_s(x,f(x))\wt{s}\\
		& \approx -\Omega^2 \wt{x} + g_s(x^*, f(x^*)) \wt{s},
	\end{split}
	\label{}
\end{equation}
where $g_s$ denotes the partial derivative of $g$ with respect to $s$ and
the consistency condition~\eqref{eqn:gconsistency} is applied.
{We then have}
\begin{equation}
	\begin{split}
	g(x,s)-s &= (g(x,f(x)+\wt{s})-f(x)) - (s - f(x)) \\
	&\approx ( g_s(x,f(x)) - 1 ) \wt{s} \\
	&\approx ( g_s(x^*, f(x^*)) - 1 ) \wt{s}.
	\end{split}
	\label{}
\end{equation}
In accord with notations used in later discussions, let us {denote}
\begin{equation}
	\mc{L} = g_s(x^*, f(x^*)),\quad \mc{K} = 1 - g_s(x^*, f(x^*)),
	\label{}
\end{equation}
with which the linearized system of Eq.~\eqref{eqn:TRBOMD_ODE1} becomes
\begin{equation}
  \frac{\rd^2}{\rd t^2}
  \begin{pmatrix}
    \wt{x}	\\
    \wt{s}
	\end{pmatrix}
	=
	\begin{pmatrix}
		-\Omega^2 & \mc{L}\\
		f'(x^*)\Omega^2 & -f'(x^*) \mc{L} - \omega^2\mc{K}
	\end{pmatrix}
	\begin{pmatrix}
		\wt{x}\\
		\wt{s}
	\end{pmatrix} 	
    := A
    \begin{pmatrix}
		\wt{x}\\
		\wt{s}
	\end{pmatrix}.
\label{eqn:linear_ODE}
\end{equation}
Note that when the force is computed
accurately, \ie
\begin{equation}
	g(x,s)=f(x), \quad \forall s,
	\label{}
\end{equation}
we have
\begin{equation}
	\mc{L}=0,\quad \mc{K}=1,
	\label{}
\end{equation}
meaning that the motion of $\wt{x}$ is decoupled from that of $\wt{s}$, and
$\wt{x}$ follows the exact harmonic motion in the linear
response regime with the accurate frequency $\Omega$.
When the force is computed
inaccurately, $\wt{x}$ is coupled with $\wt{s}$ in
Eq.~\eqref{eqn:linear_ODE}.
Actually, we can solve
\eqref{eqn:linear_ODE} analytically
and the eigenvalues of $A$ are
\begin{equation}
	\begin{pmatrix}
		\lambda_{\wt{\Omega}} \\
		\lambda_{\wt{\omega}}
	\end{pmatrix}
	=
	\begin{pmatrix}
		\frac12\left( \sqrt{(\mc{L}f'(x^*) + \mc{K}\omega^2 + \Omega^2)^2 -
		4\mc{K}\omega^2 \Omega^2} - \mc{L}f'(x^*) - \mc{K}\omega^2 -
		\Omega^2\right)\\
		\frac12\left( -\sqrt{(\mc{L}f'(x^*) + \mc{K}\omega^2 + \Omega^2)^2 -
		4\mc{K}\omega^2 \Omega^2} - \mc{L}f'(x^*) - \mc{K}\omega^2 -
		\Omega^2\right)
	\end{pmatrix}.
	\label{eqn:eig_linear_ODE}
\end{equation}
Then the frequencies of the normal modes {of the ODE} are
$\wt{\Omega}=\sqrt{-\lambda_{\wt{\Omega}}}$
and $\wt{\omega}=\sqrt{-\lambda_{\wt{\omega}}}$ respectively.
Assume
$\omega^2\gg \Omega^2$ and expand the solution to the order of
$\Or(1/\omega^2)$, we have
\begin{equation}
	\wt{\Omega} = \Omega \left( 1 - \frac{f'(x^*)}{2\omega^2}
	\mc{L}\mc{K}^{-1} \right) + \Or(1/\omega^4).
	\label{eqn:OmegaPerturbODE}
\end{equation}
Similarly the frequency for the other normal mode which is dominated by
the motion of $\wt{s}$ is
\begin{equation}
	\wt{\omega} = \sqrt{\mc{K}}\omega \left( 1 + \frac{f'(x^*)}{2\omega^2}
	\mc{L}\mc{K}^{-1}  \right) + \Or(1/\omega^3).
	\label{eqn:ficomegaPerturbODE}
\end{equation}
It is found that one of the normal mode of Eq.~\eqref{eqn:linear_ODE}
has frequency $\wt{\Omega}\approx \Omega$.  We can therefore measure the
accuracy of Eq.~\eqref{eqn:TRBOMD_ODE} using the relative error  between
$\wt{\Omega}$ and $\Omega$. Furthermore, if the
dynamics~\eqref{eqn:TRBOMD_ODE} is stable in the linear response regime,
it is necessary to have $\mc{K}>0$.

From Eq.~\eqref{eqn:OmegaPerturbODE} we conclude that if the time
reversible numerical scheme~\eqref{eqn:TRBOMD_ODE} is used for
simulating the ODE~\eqref{eqn:ODE} and if we neglect the error due to
the Verlet scheme, the error introduced in computing the frequency
$\Omega$ is proportional to $\omega^{-2}$.  This seems to indicate
that very large $\omega$ (\ie very small time step $\Delta t$) might
be needed to obtain accurate results. Fortunately the $\omega^{-2}$ term in
Eq.~\eqref{eqn:OmegaPerturbODE} has the prefactor $f'(x^{*})
\mc{L}\mc{K}^{-1}$. Eq.~\eqref{eqn:flinear} shows that
$f'(x^{*})\approx -\Omega^2$, which is small compared to $\omega^{2}$.
If $g_{s}(x^*,f(x^*))$ is small, then $\mc{K}\approx 1$, and the
accuracy of $\wt{\Omega}$ is determined by $\mc{L}$ or
$g_{s}(x^*,f(x^*))$, which indicates the sensitivity of the computed
force with respect to the initial guess, or the accuracy of the
iterative procedure for computing the force. If a ``good'' iterative
procedure is used, $g_{s}(x^*,f(x^*))$ will be small. Therefore the
presence of the term $\mc{L}$ allows one to obtain relatively
accurate approximation to the frequency $\Omega$ without using a large
$\omega$.  The same behavior can be observed when using TRBOMD to
approximate BOMD (\videpost).

Finally, we remark that even though Eq.~\eqref{eqn:ODE} is a much
simplified system, it will be seen below that for BOMD with $M$ atoms
and $N$ interacting electrons, the analysis in the linear response
regime follows the same line, and the result for the frequency is
similar to Eq.~\eqref{eqn:OmegaPerturbODE}.

\section{Time reversible Born-Oppenheimer molecular
dynamics}\label{sec:trbomd}

Consider a system with $M$ atoms and $N$ electrons.
The position of the atoms at time $t$ is denoted by
$\vR(t)=(R_{1}(t),\ldots, R_M(t))^T$. In BOMD,
the motion of atoms follows Newton's law
\begin{equation}
	m \ddot{R}_{I}(t) = f_{I}(\vR(t)) = -\frac{\partial E(\vR(t))}{\partial R_{I}},
	\label{}
\end{equation}
where $E(\vR(t))$ is the total energy of the system at the atomic
configuration $\vR(t)$.
In KSDFT, the total energy is expressed as a
functional of a set of Kohn-Sham orbitals $\{\psi_i(x)\}_{i=1}^{N}$.
To illustrate the idea with minimal technicality, let us consider for
the moment a system of $N$ electrons at zero temperature.
The energy functional in KSDFT takes the form
\begin{equation}
	\begin{split}
  E(\{\psi_{i}(x)\}_{i=1}^{N};\vR)
  =& \frac12 \sum_{i=1}^{N}\int \abs{\nabla \psi_i(x)}^2\ud x
	+ \int  \rho(x) V_{\ion}(x;\vR)\ud x + E_{\hxc}[\rho],\\
	\rho(x)=& \sum_{i=1}^{N} \abs{\psi_{i}(x)}^{2}.
	\end{split}
  \label{eqn:Esimplify}
\end{equation}
The first term in the energy functional is
the kinetic energy of the electrons. The second term contains the electron-ion
interaction energy. The ion-ion interaction energy usually takes the
form $\sum_{I< J} \frac{Z_I Z_J}{\abs{R_I-R_J}}$ where $Z_{I}$ is the
charge for the nucleus $I$.  The ion-ion interaction energy does not depend
on the electron density $\rho$. To simplify the notation,
we include the ion-ion interaction energy in the $V_{\ion}$ term as a
constant shift that is independent of the $x$ variable. The third
term does not explicitly depend on the atomic configuration $\vR$, and
is a nonlinear functional of the electron density $\rho$. It represents
the Hartree part of electron-electron interaction energy (h), and the
exchange-correlation energy (xc) characterizing  many body effects.
The energy $E(\vR)$ as a function of
atomic positions is given by the following minimization problem
\begin{equation}
  \begin{split}
    & E(\vR) = \min_{\{\psi_{i}(x)\}_{i=1}^{N}}
    E(\{\psi_{i}(x)\}_{i=1}^{N};\vR), \\
    &\text{s.t.} \quad \int \psi_{i}^{\dagger}(x) \psi_{j}(x)\ud x = \delta_{ij},
    \quad i,j=1,\ldots,N.
  \end{split}
  \label{eqn:Eminpsi}
\end{equation}
We denote by $\{\psi_{i}(x;\vR)\}_{i=1}^{N}$ the (local) minimizer, and
$\rho^{*}(x;\vR)=\sum_{i=1}^{N}\abs{\psi_{i}(x;\vR)}^2$ the
converged electron density corresponding to the minimizer
(here we assume that the minimizing electron density is unique).
Then the force acting on the atom $I$ is
\begin{equation}
	f_{I}( \vR; \rho^{*}(x;\vR) ) =
	-\frac{\partial E(\vR)}{\partial R_{I}}
	= -\int \rho^{*}(x;\vR)
	\frac{\partial V_{\ion}(x;\vR)}{\partial R_{I}} \ud x.
	\label{eqn:HellmannFeynman}
\end{equation}
In physics literature the force formula in
Eq.~\eqref{eqn:HellmannFeynman} is referred to as the Hellmann-Feynman
force. The validity of the Hellmann-Feynman formula relies on
that the electron density $\rho^{*}(x;\vR)$ corresponds to the
minimizers of the Kohn-Sham energy functional.  Since
$E_{\hxc}[\rho]$ is a nonlinear functional of $\rho$, the electron
density $\rho$ is usually determined through the self-consistent field
(SCF) iteration as follows.

Starting from an inaccurate input electron density
$\rho^{\mathrm{in}}$, one first computes
the output electron density by solving the lowest $N$ eigenfunctions of
the problem
\begin{equation}
	\left(-\frac12 \Delta_x + \mc{V}(x;\vR,\rho^{\mathrm{in}})\right)\psi_{i} =
	\varepsilon_{i} \psi_{i}
	\label{}
\end{equation}
with
\begin{equation}
	\mc{V}(x;\vR,\rho) = V_{\ion}(x;\vR) + \frac{\delta
	E_{\hxc}[\rho]}{\delta \rho}(x),
	\label{}
\end{equation}
and the output electron density $\rho^{\mathrm{out}}$ is defined by
\begin{equation}
	\rho^{\mathrm{out}}(x) := F[\rho^{\mathrm{in}}](x) =
	\sum_{i=1}^{N} \abs{\psi_{i}(x)}^2.
	\label{}
\end{equation}
Here the operator $F$ is called the Kohn-Sham map.
$\rho^{\mathrm{out}}$ can be used directly as the input electron density
$\rho^{\mathrm{in}}$ in the next iteration.  This is called the \textit{fixed
point iteration}. Unfortunately, in most electronic structure
calculations, the fixed point iteration does not converge even when
$\rho^{\mathrm{in}}$ is very close to the true electron density
$\rho^{*}$.  The fixed point iteration can be improved by the simple
mixing method, which takes the linear combination of the electron
density
\begin{equation}
	\alpha\rho^{\mathrm{out}} + (1-\alpha)\rho^{\mathrm{in}}
	\label{}
\end{equation}
as the input density for the next iteration with $0<\alpha\le
1$. Simple mixing can greatly improve the convergence properties of
the SCF iteration over the fixed point iteration,
but the convergence rate can still be slow in practice.
There are more complicated SCF iteration schemes
such as Anderson mixing scheme~\cite{Anderson1965}, Pulay mixing
scheme~\cite{Pulay1980} and Broyden mixing scheme~\cite{Johnson1988}.
Furthermore, preconditioners can be applied to the SCF iteration to
enhance convergence properties such as the Kerker
preconditioner~\cite{Kerker1981}.  More detailed discussion on
convergence properties of these SCF schemes can be found
in~\cite{LinYang2013}. In the following discussions, we denote by
$\rhoscf(x;\vR, \rho)$ the final electron density after the SCF
iteration starting from an initial guess $\rho$. We assume that
$\rhoscf$ satisfies the consistency condition
\begin{equation}
	\rhoscf(x;\vR,\rho^{*}(\cdot;\vR)) = \rho^{*}(x;\vR).
	\label{eqn:rhoconsistency}
\end{equation}
If a non-convergent SCF iteration procedure is used, $\rhoscf(x;\vR,
\rho)$ might deviate from $\rho^{*}(x;\vR)$. Such deviation introduces
error in the force, and the error can accumulate in the long time
molecular dynamics simulation, and lead to inaccurate results in
computing the statistical and dynamical properties of the systems.

The map $\rhoscf$ is usually highly nonlinear, 
which makes it
difficult to correct the error in the force.  The TRBOMD scheme
avoids the direct
correction for the inaccurate $\rhoscf$, but allows the initial guess to
dynamically evolve together with the motion of the atoms.  We
denote by $\rho(x,t)$ the initial guess for the
SCF iteration at time $t$.  When $\rho(\cdot,t)$ is used as an argument, we
also write
$\rhoscf(x;\vR(t),\rho(t)):=\rhoscf(x;\vR(t),\rho(\cdot,t))$.
The Hellmann-Feynman
formula~\eqref{eqn:HellmannFeynman} is used to
compute the force at the electron density $\rhoscf(x;\vR(t),\rho(t))$ even
though $\rho^{*}(x;\vR(t))$ is not available.  Thus,
the equation of motion in TRBOMD reads
\begin{equation}
	\begin{split}
	m \ddot{R}_{I}(t) &= f_{I}( \vR(t); \rhoscf(x; \vR(t),\rho(t) )) = -\int
	\rhoscf(x; \vR(t),\rho(t) )	\frac{\partial
	V_{\ion}(x;\vR(t))}{\partial R_I} \ud x, \\
	\ddot{\rho}(x,t) &= \omega^{2}( \rhoscf(x; \vR(t),\rho(t) ) -
	\rho(x,t) ).
	\end{split}
	\label{eqn:TRBOMD}
\end{equation}
It is clear that TRBOMD is time reversible.  The discretized TRBOMD is
still time reversible if the numerical scheme is time reversible.
For instance, if the Verlet scheme is used, the
discretized equation of motion becomes
\begin{equation}
\begin{split}
	R_{I}(t_{k+1}) &= 2 R_{I}(t_{k}) - R_{I}(t_{k-1}) -
	\frac{\Delta t^2}{m}f_{I}( \vR(t_{k}); \rhoscf(x; \vR(t_{k}),\rho(t_{k})
	),\\
	\rho(x,t_{k+1}) &= 2 \rho(x,t_{k}) - \rho(x,t_{k-1}) +
	\Delta t^2\omega^{2}( \rhoscf(x; \vR(t_{k}),\rho(t_{k}) ) -
	\rho(x,t_{k}) ),
\end{split}\label{eqn:Verlet}
\end{equation}
which is evidently time reversible.  The artificial frequency $\omega$ controls the frequency
of the fictitious dynamics of $\rho(x,t)$ and is generally chosen to be
larger than the frequency of motion of the atoms.
The numerical stability of the Verlet algorithm requires that the
dimensionless quantity $\kappa:=(\omega\Delta t)^2$ to be
small~\cite{McLachlanAtela1992}.  When $\kappa$ is fixed, $\omega$
controls the stiffness, or equivalently the time step $\Delta
t=\frac{\sqrt{\kappa}}{\omega} $ for the equation of
motion~\eqref{eqn:Verlet}.

Let us mention that TRBOMD is closely related to CPMD.
In CPMD, the equation of motion is given by
\begin{equation}
\begin{split}
  m \ddot{R}_I(t) & = f_I(\bd{R}(t), \rho(t)) = - \int \rho(t)
  \frac{\partial V_{\text{ion}}(x; \bd{R}(t))}{\partial R_I} \ud x, \\
  \mu \ddot{\psi}_i(t) & = - \frac{\delta E(\bd{R}(t),
    \{\psi_i(t)\})}{\delta \psi_i^{\dagger}} + \sum_{j} \psi_j(t)
  \Lambda_{ji}(t),
\end{split}\label{eqn:CPMD}
\end{equation}
where $\mu$ is the fictitious electron mass for the fake electron
dynamics in CPMD, and $\Lambda$'s are the Lagrange multipliers determined
so that $\{\psi_i(t)\}$ is an orthonormal set of functions for any
time. The CPMD scheme~\eqref{eqn:CPMD} can be viewed as the equation of
motion with an extended Lagrangian
\begin{equation}
  \mc{L}_{\text{CP}}\bigl(\bd{R}, \dot{\bd{R}}, \{\psi_i\}, \{\dot{\psi}_i\}\bigr)
  = \sum_I \frac{m}{2} \lvert \dot{R}_I\rvert^2 + \sum_i \frac{\mu}{2}
  \int \lvert \dot{\psi}_i\rvert^2 - E(\bd{R}, \{\psi_i\}),
\end{equation}
which contains both ionic and electronic degrees of
freedom. Therefore, CPMD is a Hamiltonian dynamics and thus time
reversible.

Note that the frequency of the evolution equation for $\{\psi_{i}\}$ in CPMD
is adjusted by the fictitious mass parameter $\mu$. Comparing with
TRBOMD, the parameter $\mu$ plays a similar role as $\omega^{-2}$ which
controls the frequency of the fictitious dynamics of the initial density
guess in SCF iteration. This connection will be made more explicit in
the sequel.

We remark that the papers~\cite{Niklasson2008,Niklasson2009} made a
further step in viewing TRBOMD by an extended Lagrangian approach in a
vanishing mass limit. However, unless very specific and restrictive
form of the error due to non-convergent SCF iterations is assumed, the
equation of motion in TRBOMD does not have an associated Lagrangian in
general.  The connection remains formal, and hence we will not further
explore here.

\section{Analysis of TRBOMD in the linear response regime}\label{sec:lr}

In this section we consider Eq.~\eqref{eqn:TRBOMD} in the linear
response regime, in which each atom $I$ oscillates
around its equilibrium position $R_{I}^{*}$.  The displacement of the atomic configuration $\vR$ from the
equilibrium position is denoted by $\wt{\vR}(t):=\vR(t)-\vR^{*}$, and the deviation of the
electron density from the converged density is denoted by
$\wt{\rho}(x,t):=\rho(x,t) - \rho^{*}(x;\vR(t))$.
Both $\wt{\vR}(t)$ and $\wt{\rho}(x,t)$ are small quantities in the linear response
regime, and contain the same information as $\vR(t)$ and $\rho(x,t)$.
Using $\wt{\vR}(t)$ and $\wt{\rho}(x,t)$ as the new variables and
noting the chain rule due to the
$\vR$-dependence in $\rho^{*}(x;\vR(t))$, the
equation of motion in TRBOMD becomes
\begin{equation}
	\begin{split}
		m \ddot{\wt{R}}_{I}(t) & = -\int
	\rhoscf(x; \vR(t),\rho(t) )	\frac{\partial
	V_{\ion}(x;\vR(t))}{\partial R_I} \ud x, \\
	\ddot{\wt{\rho}}(x,t) & = \omega^{2}( \rhoscf(x; \vR(t),\rho(t) ) -
	\rho(x,t) ) - \sum_{I=1}^{M}\frac{\partial
	\rho^{*}(x;\vR(t))}{\partial R_{I}}
	\ddot{\wt{R}}_{I}(t) \\
	&\qquad - \sum_{I,J=1}^{M} \dot{\wt{R}}_{I}(t)\dot{\wt{R}}_{J}(t)
	\frac{\partial^2 \rho^{*}(x;\vR(t))}{\partial R_I \partial
	R_{J}}.
	\end{split}
	\label{eqn:residTRMD}
\end{equation}
To simplify notation from now on we suppress the $t$-dependence in all
variables, and Eq.~\eqref{eqn:residTRMD} becomes
\begin{subequations}\label{eqn:residTRMD2}
\begin{align}
      m \ddot{\wt{R}}_{I} & = -\int \rhoscf(x; \vR,\rho )
      \frac{\partial V_{\ion}(x;\vR)}{\partial R_I} \ud x, \label{eqn:residTRMD2-a}\\
      \ddot{\wt{\rho}}(x) & = \omega^{2}( \rhoscf(x; \vR,\rho ) -
      \rho(x) ) - \sum_{I=1}^{M}\frac{\partial \rho^{*}}{\partial
        R_{I}}(x;\vR)
      \ddot{\wt{R}}_{I}  - \sum_{I,J=1}^{M} \dot{\wt{R}}_{I}\dot{\wt{R}}_{J}
      \frac{\partial^2 \rho^{*}}{\partial R_I \partial R_{J}}(x;\vR). \label{eqn:residTRMD2-b}
\end{align}
\end{subequations}
In the linear response regime, we expand Eq.~\eqref{eqn:residTRMD2} and
only keep terms that are linear with respect to $\wt{\vR}$ and $\wt{\rho}$.
All the higher order terms, including all the cross products of $\wt{R}_{I}$,
$\dot{\wt{R}}_{I}$, and $\wt{\rho}$ will be dropped.
First we linearize the force on atom $I$
with respect to $\wt{\rho}$ as
\begin{equation}
	\begin{split}
  &f_{I}( \vR; \rhoscf(x; \vR,\rho )) \\
  =& -\int
	\rhoscf(x; \vR,\rho )	\frac{\partial
	V_{\ion}(x;\vR)}{\partial R_I} \ud x\\
	=& -\int \rho^{*}(x;\vR)  \frac{\partial
	V_{\ion}(x;\vR)}{\partial R_I} \ud x
	- \int \left( \rhoscf(x; \vR,\rho^{*}(\vR) + \wt{\rho} ) - \rho^{*}(x;\vR) \right) \frac{\partial
	V_{\ion}(x;\vR)}{\partial R_I} \ud x \\
	\approx&
	-\int \rho^{*}(x;\vR)  \frac{\partial
	V_{\ion}(x;\vR)}{\partial R_I} \ud x
	- \int \frac{\delta\rhoscf}{\delta\rho}(x,y;\vR) \wt{\rho}(y) \frac{\partial
	V_{\ion}(x;\vR)}{\partial R_I} \ud x \ud y.
	\end{split}
	\label{eqn:fapprox1}
\end{equation}
Next we linearize with respect to $\wt{\vR}$, we have
\begin{equation}
	\int \rho^{*}(x;\vR)  \frac{\partial
	V_{\ion}(x;\vR)}{\partial R_I} \ud x  \approx
-m\sum_{I,J=1}^{M}
	\mc{D}_{IJ} \wt{R}_{J}.
	\label{eqn:LinearForce}
\end{equation}
Here the matrix $\{\mc{D}_{IJ}\}$ is the dynamical matrix for the atoms. For
the last term in Eq.~\eqref{eqn:fapprox1} we have
\begin{equation}
	\begin{split}
	&\int \frac{\delta\rhoscf}{\delta\rho}(x,y;\vR) \wt{\rho}(y) \frac{\partial
	V_{\ion}(x;\vR)}{\partial R_I} \ud x \ud y \\
	\approx &
	\int \frac{\delta\rhoscf}{\delta\rho}(x,y;\vR^{*}) \wt{\rho}(y) \frac{\partial
	V_{\ion}(x;\vR^{*})}{\partial R_I} \ud x \ud y \\
	:= & -m \mc{L}_{I}[\wt{\rho}].
	\end{split}
	\label{eqn:linear1}
\end{equation}
The last equation in Eq.~\eqref{eqn:linear1} defines a linear functional
$\mc{L}_{I}$, with $\frac{\delta\rhoscf}{\delta \rho}(x,y;\vR^{*})$ and
$\frac{\partial V_{\ion}(x;\vR^{*})}{\partial R_I}$ evaluated at the
fixed
equilibrium point $\vR^{*}$.

In the linear response regime, the operator
$\frac{\delta\rhoscf}{\delta\rho}(x,y;\vR^{*})$ carries all the
information of the SCF iteration scheme. Let us now derive the explicit
form of $\frac{\delta\rhoscf}{\delta\rho}(x,y;\vR^{*})$ for
the $k$-step simple mixing scheme with mixing parameter (step length)
$\alpha$ ($0<\alpha\le 1$).  If $k=1$, the simple mixing scheme
reads
\begin{equation}
	\rhoscf( x; \vR, \rho^{*}(\vR) + \wt{\rho} ) =
	\alpha F[\rho^{*}(\vR)+\wt{\rho}] + (1-\alpha)
	(\rho^{*}(\vR)+\wt{\rho}),
	\label{}
\end{equation}
so
\begin{equation}
	\frac{\delta\rhoscf}{\delta \rho}(x,y;\vR^{*})
	= \delta(x-y) - \alpha \left(\delta(x-y) - \frac{\delta F}{\delta
	\rho}(x,y)\right).
	\label{eqn:simpleone}
\end{equation}
Here $\delta(x)$ is the Dirac $\delta$-function, and the operator
$\left(\delta(x-y) - \frac{\delta F}{\delta \rho}(x,y)\right):= \varepsilon(x,y)$  is usually refereed to as the
\textit{dielectric operator}~\cite{Adler1962,Wiser1963}.
To simplify the notation we would not distinguish the kernel of an
integral operator from the integral operator itself. For example
$\varepsilon(x,y)$ is denoted by $\varepsilon$.
Neither will we distinguish integral operators defined on continuous
space from the corresponding finite dimensional matrices obtained from
certain numerical discretization. This slight abuse of notation allows
us to simply denote $f(x)=\int A(x,y) g(y)\ud y$ by $f=Ag$ as
a matrix-vector multiplication, and to denote the composition of kernels
of integral operators $C(x,y)=\int dz A(x,z) B(z,y)$ by $C=AB$ as
a matrix-matrix multiplication.  Using such notations,
Eq.~\eqref{eqn:simpleone} can be written in a more compact form
\begin{equation}
	\frac{\delta\rhoscf}{\delta\rho} = I - \alpha \varepsilon.
	\label{}
\end{equation}
Similarly for the $k$-step simple mixing method, we have
\begin{equation}
	\frac{\delta\rhoscf}{\delta\rho} = (1-\alpha
	\varepsilon)^{k}.
	\label{eqn:simplek}
\end{equation}
In general the dielectric operator is diagonalizable and all eigenvalues of
$\varepsilon$ are real.  Therefore the linear response operator
$\frac{\delta\rhoscf}{\delta\rho}$ for the $k$-th step simple
mixing method is also diagonalizable with real eigenvalues.

From Eq.~\eqref{eqn:residTRMD2-b} we have
\begin{equation}
	\begin{split}
	& \rhoscf(x; \vR,\rho ) - \rho(x)	 \\
	=& \left( \rhoscf(x; \vR,\wt{\rho} + \rho^{*}(\vR) ) - \rho^{*}(x;\vR) \right)
	- \left( \rho(x) - \rho^{*}(x;\vR \right)) \\
	\approx& \int \frac{\delta \rhoscf}{\delta \rho}(x,y;\vR)
	\wt{\rho}(y) \ud y - \wt{\rho}(x)\\
	\approx& \int \frac{\delta \rhoscf}{\delta \rho}(x,y;\vR^*)
	\wt{\rho}(y) \ud y - \wt{\rho}(x)\\
	:=& - \int \mc{K}(x,y) \wt{\rho}(y) \ud y.
	\end{split}
	\label{eqn:linear2}
\end{equation}
Here we have used the consistency condition~\eqref{eqn:rhoconsistency}.
The last line of Eq.~\eqref{eqn:linear2} defines a kernel
\begin{equation}
	\mc{K}(x,y) = \delta(x-y) -  \frac{\delta\rhoscf}{\delta
	\rho}(x,y;\vR^{*}),
	\label{eqn:Kkernel}
\end{equation}
which is an important quantity for the stability of TRBOMD as
will be seen later. Using Eqs.~\eqref{eqn:linear1} and \eqref{eqn:linear2},
the equation of motion \eqref{eqn:residTRMD2} can be written in the linear
response regime as
\begin{equation}
	\begin{split}
		\ddot{\wt{R}}_{I}      &= -\sum_{J=1}^{M} \mc{D}_{IJ} \wt{R}_{J} +
	\mc{L}_{I}[\wt{\rho}], \\
	\ddot{\wt{\rho}}(x) &= -\omega^{2}\int \mc{K}(x,y) \wt{\rho}(y) \ud y -
	\sum_{I=1}^{M} \frac{\partial \rho^{*}}{\partial R_{I}}(x;\vR^*)
	\left(-\sum_{J=1}^{M} \mc{D}_{IJ} \wt{R}_{J} +
	\mc{L}_{I}[\wt{\rho}]\right).
	\end{split}
	\label{eqn:linearTRBOMD1}
\end{equation}
Define
\begin{equation}
	\vL = (\mc{L}_{1}, \cdots, \mc{L}_{M})^{T},
	\label{}
\end{equation}
then Eq.~\eqref{eqn:linearTRBOMD1} can be rewritten
in a more compact form as
\begin{subequations}\label{eqn:linearTRBOMD2}
	\begin{align}
		\ddot{\wt{\vR}}      &= -\mc{D} \wt{\vR} + \vL[\wt{\rho}], \label{eqn:linearTRBOMD2-a}\\
	\ddot{\wt{\rho}}(x) &= -\omega^{2}\int \mc{K}(x,y) \wt{\rho}(y) \ud y -
	\left(\frac{\partial \rho^{*}}{\partial \vR}(x;\vR^*)\right)^T
	\left(- \mc{D} \wt{\vR} + \vL[\wt{\rho}]\right).\label{eqn:linearTRBOMD2-b}
	\end{align}
\end{subequations}

Now if the self-consistent iteration is performed accurately regardless
of the initial guess, \ie
\begin{equation}
	\rhoscf(x;\vR,\rho) = \rho^{*}(x;\vR), \quad \forall \rho,
	\label{}
\end{equation}
which implies
\begin{equation}
	\frac{\delta \rhoscf}{\delta \rho}(x,y;\vR^*) = 0,\quad  \vL = \mathbf{0}, \quad \mc{K}(x,y) = \delta(x-y).
	\label{}
\end{equation}
The linearized equation of
motion~\eqref{eqn:linearTRBOMD2} becomes
\begin{subequations}\label{eqn:linearTRBOMD3}
	\begin{align}
		\ddot{\wt{\vR}}      &= -\mc{D} \wt{\vR}, \label{eqn:linearTRBOMD3-a}\\
		\ddot{\wt{\rho}}(x) &= -\omega^{2}\wt{\rho}(x) +
	\left(\frac{\partial \rho^{*}}{\partial \vR}(x;\vR^*)\right)^T
	\mc{D} \wt{\vR}.\label{eqn:linearTRBOMD3-b}
	\end{align}
\end{subequations}
Therefore in the case of accurate SCF iteration, according to Eq.~\eqref{eqn:linearTRBOMD3-a},
the equation of motion
of atoms follows the accurate linearized equation, and is decoupled
from the fictitious dynamics of $\wt{\rho}$.
The normal modes of the
equation of motion of atoms can be obtained by diagonalizing the
dynamical matrix $\mc{D}$ as
\begin{equation}
	\mc{D} \vv_{l} = \Omega_{l}^2 \vv_{l},\quad l=1,\ldots,M.
	\label{eqn:phonon1}
\end{equation}
The frequencies $\{\Omega_{l}\}$ ($\Omega_{l}>0$) are known as
\textit{phonon frequencies}. When the SCF iterations are performed inaccurately,
it is meaningless to
assess the accuracy of the approximate
dynamics~\eqref{eqn:linearTRBOMD2} by direct investigation of the
trajectories $\wt{R}(t)$,
since small difference
in the phonon frequency can cause large error in the phase of the periodic motion
$\wt{R}(t)$ over long time. However, it is possible to compute the
approximate phonon frequencies $\{\wt{\Omega}_{l}\}$ from
Eq.~\eqref{eqn:linearTRBOMD2}, and measure the accuracy of TRBOMD
in the linearized regime from the relative error
\begin{equation}
	\err_{l} = \frac{\wt{\Omega}_{l}-\Omega_{l}}{\Omega_{l}}.
	\label{eqn:errorphonon}
\end{equation}

The operator $\mc{K}(x,y)$ in Eq.~\eqref{eqn:Kkernel} is directly related to the stability of the
dynamics.  Eq.~\eqref{eqn:linearTRBOMD2-b} also suggests that in the
linear response regime, the spectrum of $\mc{K}(x,y)$ must be on the
real line, which requires that the matrix
$\frac{\delta\rhoscf}{\delta\rho}(x,y;\vR^{*})$ be diagonalizable
with real eigenvalues. This has been shown for the simple mixing
scheme.  However, we remark that the condition that all eigenvalues of
$\mc{K}(x,y)$ are real may not hold for general preconditioners or for
more complicated SCF iterations (for instance,
Anderson mixing).  This is one important restriction of the linear
response analysis. Of course, this may not be a restriction for
practical TRBOMD simulation for real systems. We will leave further
understanding of this to future works.

Let us now assume that all eigenvalues of $\mc{K}$ are real.  The lower bound
of the spectrum of $\mc{K}$, denoted by $\lambda_{\min}(\mc{K})$, should
satisfy
\begin{equation}
	\lambda_{\min}(\mc{K})>0.
	\label{eqn:stability}
\end{equation}
Eq.~\eqref{eqn:stability} is a necessary condition for TRBOMD
to be stable, which will be referred to as the \textit{stability
condition} in the following.  Furthermore, $\omega$ should be
chosen large enough in order to avoid resonance between the motion
of $\wt{\vR}$ and $\wt{\rho}$.
Therefore the \textit{adiabatic condition}
\begin{equation}
	\omega^2 \gg \frac{\lambda_{\max}(\mc{D})}{\lambda_{\min}(\mc{K})} =
	\frac{\max_{l} \Omega_{l}^2}{\lambda_{\min}(\mc{K})}
	\label{eqn:stable1}
\end{equation}
should also be satisfied.  Due to Eq.~\eqref{eqn:stable1}, we may assume
$\epsilon=1/\omega^2$ is a small number, and expand $\Omega_{l}$
in the perturbation series of $\epsilon$ to quantify the error in
the linear response regime.  Following the derivation in the appendix, we have
\begin{equation}
	\wt{\Omega}_{l} = \Omega_{l}\left( 1 - \frac{1}{2\omega^2}
	\vv_{l}^T \vL\left[ \mc{K}^{-1}\left[ \left(\frac{\partial
	\rho^{*}}{\partial \vR}\right)^T \vv_{l} \right] \right]
	\right) + \Or(1/\omega^4),
	\label{eqn:OmegaPerturb}
\end{equation}
where $\mc{K}^{-1}$ is the inverse operator of $\mc{K}$ ($\mc{K}$ is invertible due to the stability condition).
Since $\omega=\sqrt{\kappa}/\Delta t$, Eq.~\eqref{eqn:OmegaPerturb} suggests
that the accuracy of TRBOMD in the linear response regime is
$(\Delta t)^2$, with preconstant mainly determined by $\vL$, \ie the
accuracy of the SCF iteration.


Let us compare TRBOMD with CPMD.  It is well
known that CPMD accurately approximates the results of BOMD, provided
that the electronic and ionic degrees of freedom remain adiabatically
separated as well as the electrons stay close to the Born-Oppenheimer
surface \cite{PastoreSmargiassiBuda1991,BornemannSchutte1998}. More
specifically, the fictitious electron mass should be chosen so that
the lowest electronic frequency is well above ionic frequencies
\begin{equation}\label{eq:CPMDadiabatic}
  \mu \ll  \frac{E_{\text{gap}}}{\max_l \Omega_l^2},
\end{equation}
where $E_{\text{gap}}$ is the spectral gap (between highest occupied
and lowest unoccupied states) of the system, and recall that
$\Omega_l$ is the vibration frequency of the lattice phonon.
For CPMD,
a similar analysis in the linear response regime as above (which we omit the derivation here) shows that
\begin{equation}\label{eq:CPMDaccuracy}
  \wt{\Omega}_l = \Omega_l (1 + \Or(\mu)),
\end{equation}
under the assumption~\eqref{eq:CPMDadiabatic}.

Note that the condition~\eqref{eq:CPMDadiabatic}
implies that CPMD  no longer works if the
system has a small gap or is even metallic. The usual work-around for
this is to add a heat bath for the electronic degrees of freedom in
CPMD~\cite{BlochlParrinello1992}, so that it maintains a fictitious
temperature for the electronic degree of freedom. Nonetheless the
adiabaticity is lost  for metallic systems and CPMD is no longer
accurate over long time simulation. In contrast, as we have discussed
previously, TRBOMD may work for both insulating and
metallic systems without any modification, provided that the SCF
iteration is accurate and no resonance occurs. This is an important
advantage of TRBOMD, which we will illustrate using numerical examples
in the next section.

When the system has a gap we can take $\mu$ sufficiently small
to satisfy the adiabatic separation condition \eqref{eq:CPMDadiabatic}.
Compare Eq.~\eqref{eq:CPMDaccuracy} with Eq.~\eqref{eqn:OmegaPerturb},
we see that $\mu$ in CPMD plays a similar role as $\omega^{-2}$ in TRBOMD. The accuracy (in the linear regime) for CPMD and TRBOMD is first order in $\mu$ and $\omega^{-2}$ respectively. At the same time, as taking a small $\mu$ or large $\omega$ increases the stiffness of the equation, the computational cost is proportional to $\mu^{-1}$ and $\omega^2$, respectively.

Let us remark that the above analysis is done in the linear response
regime. As shown in \cite{PastoreSmargiassiBuda1991,BornemannSchutte1998}, the accuracy of CPMD in general is only $\Or(\mu^{1/2})$ instead of $\Or(\mu)$ for the linear regime. Due to the close connection between these two parameters, we do not expect $\Or(\omega^{-2})$ accuracy for TRBOMD in general either. Actually, as will be discussed in Section~\ref{sec:nonlinear},
{if the deviation of atom positions from equilibrium is not so small
that we cannot linearize the nuclei motion, the error of TRBOMD in
general will be $\Or(\omega^{-1})$.}

\section{Numerical results in the linear response
regime}\label{sec:numer}

In this section we present numerical results for TRBOMD in the linear
response regime using a one dimensional (1D) model for KSDFT without
the exchange correlation functional.  The model
problem can be tuned to exhibit both metallic and insulating
features.  Such model was used before in mathematical analysis of ionization conjecture \cite{Solovej1991}.

The total energy functional in our 1D density functional theory (DFT) model is given by
\begin{equation}
  E(\{\psi_{i}(x)\}_{i=1}^{N};\vR)
	= \frac12 \sum_{i=1}^{N}\int \abs{\frac{d}{dx} \psi_i(x)}^2\ud x
	+ \frac12 \int K(x,y) (\rho(x)+m(x;\vR)) (\rho(y)+m(y;\vR)) \ud x \ud y,
  \label{eqn:E1D}
\end{equation}
with $\rho(x)= \sum_{i=1}^{N} \abs{\psi_{i}(x)}^{2}$. The associated
Hamiltonian is given by
\begin{equation}
   H(\vR) = -\frac{1}{2} \frac{d^2}{dx^2} + \int K(x,y)
	 (\rho(y)+m(y;\vR)) \ud y.
   \label{eqn:HrHF}
\end{equation}
Here $m(x;\vR)=\sum_{I=1}^{M} m_I(x-R_I)$, with the position of the
$I$-th nucleus denoted by $R_{I}$.  Each function $m_{I}(x)$ takes the
form
\begin{equation}
  m_{I}(x) = -\frac{Z_I}{\sqrt{2\pi\sigma_{I}^2}}
  e^{-\frac{x^2}{2\sigma_I^2}},
  \label{}
\end{equation}
where $Z_I$ is an integer representing the charge of the $i$-th nucleus.
This can be understood as a local pseudopotential approximation to
represent the electron-ion interaction. The second term on the right
hand side of Eq.~\eqref{eqn:E1D} represents the electron-ion,
electron-electron and ion-ion interaction energy.
The parameter $\sigma_{I}$ represents the width of the nuclei in the
pseudopotential theory.  Clearly as $\sigma_{I}\to 0$,
$m_{I}(x)\to -Z_{I}\delta(x)$ which is the
charge density for an ideal nucleus. In our numerical simulation, we set
$\sigma_{I}$ to a finite value. The corresponding $m_{I}(x)$ is
called a \textit{pseudo charge density} for the $I$-th nucleus. We
refer to the function $m(x)$ as the total pseudo-charge density of the nuclei.
The system satisfies charge
neutrality condition, \ie
\begin{equation}
  \int \rho(x) + m(x;\vR) \ud x = 0.
  \label{eqn:chargeneutral}
\end{equation}
Since $\int m_I(x) \ud x = -Z_I$, the charge neutrality
condition~\eqref{eqn:chargeneutral} implies
\begin{equation}
  \int \rho(x) \ud x = \sum_{I=1}^{M} Z_{I} = N,
  \label{}
\end{equation}
where $N$ is the total number of electrons in the system. To
simplify discussion, we omit the spin degeneracy here. The
Hellmann-Feynman force is given by
\begin{equation}
	f_I = -\int K(x,y) (\rho(y)+m(y;\vR)) \frac{\partial
	m(x;\vR)}{\partial R_I} \ud x \ud y.
	\label{eqn:force1D}
\end{equation}

Instead of using a bare Coulomb interaction, which diverges in
1D, we adopt a Yukawa kernel
\begin{equation}
  K(x,y) = \frac{2\pi e^{-\kappa \abs{x-y}}}{\kappa\epsilon_{0}},
  \label{eqn:YukawaK1}
\end{equation}
which satisfies the equation
\begin{equation}
  -\frac{d^2}{d x^2} K(x,y) + \kappa^2 K(x,y) = \frac{4\pi}{\epsilon_0} \delta(x-y).
  \label{eqn:YukawaK2}
\end{equation}
As $\kappa\to 0$, the Yukawa kernel approaches the bare Coulomb interaction
given by the Poisson equation. The parameter $\epsilon_0$ is used to
make the magnitude of the electron static contribution comparable to
that of the
kinetic energy.

\begin{figure}[h]
\centering
\subfigure[Insulator.]
{\includegraphics[width=8cm,height=6cm]{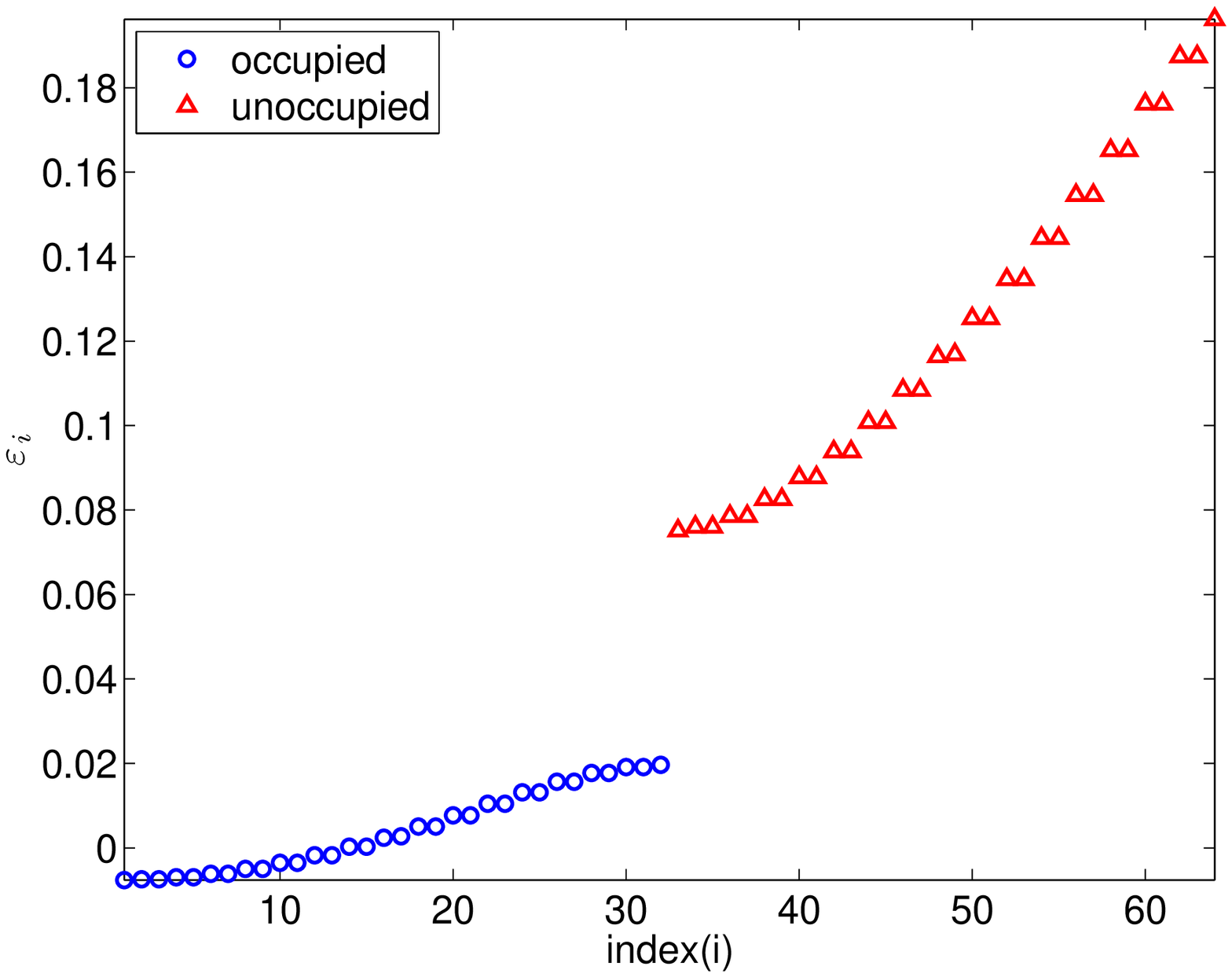}}
\subfigure[Metal.]
{\includegraphics[width=8cm,height=6cm]{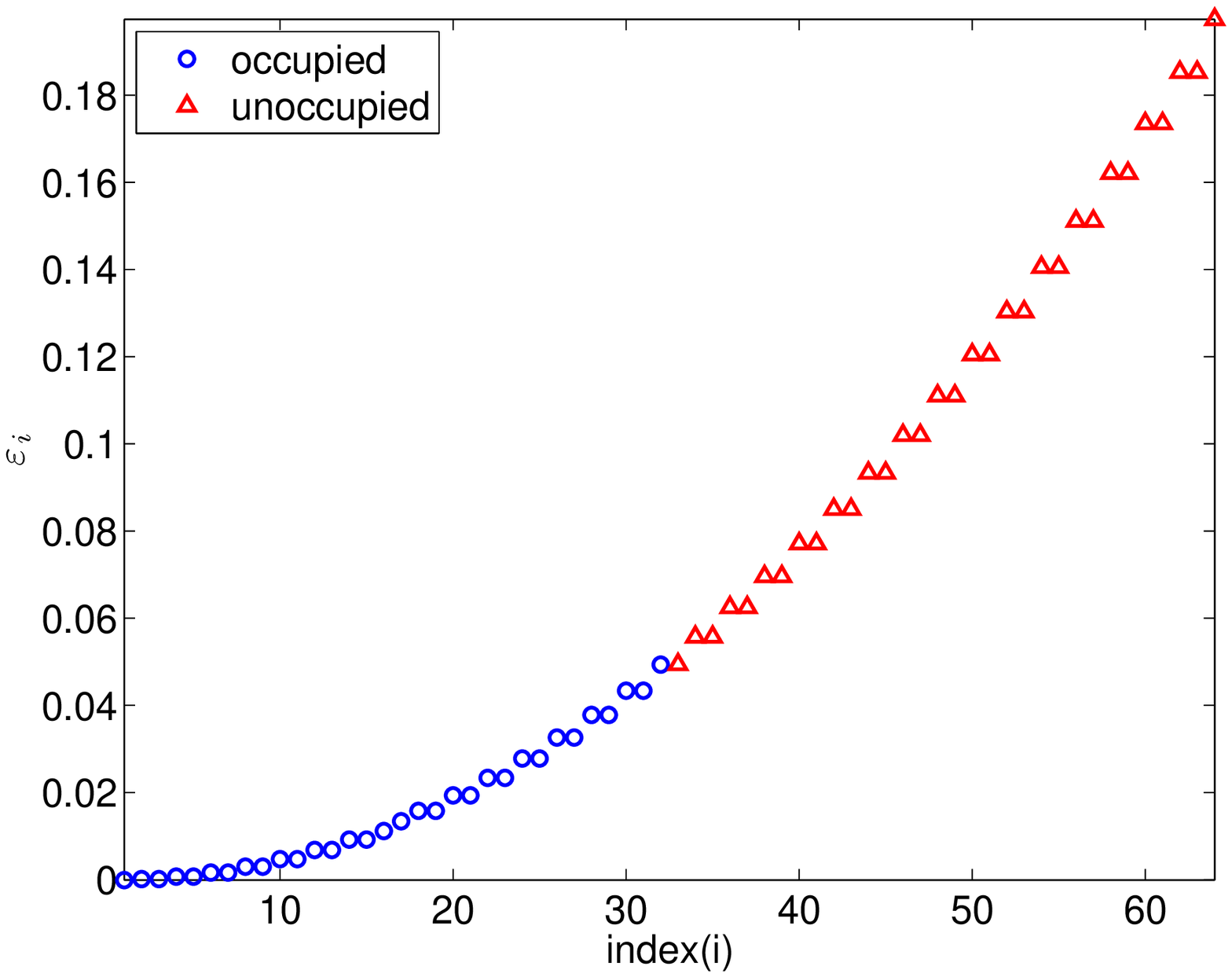}}
\caption{Spectrum for insulator and metal with $32$ atoms after $1000$ BOMD steps with converged SCF iteration.}
	\label{fig:eigH}
\end{figure}

The parameters used in the 1D DFT model are chosen as
follows.  Atomic units are used throughout the discussion unless
otherwise mentioned.  The
Yukawa parameter $\kappa=0.01$ is small enough so that the range
of the electrostatic interaction is sufficiently long, and $\epsilon_0$
is set to $10.00$.  The nuclear charge $Z_{I}$ is set to $1$ for all atoms.
Since spin is neglected, $Z_{I}=1$ implies that each atom contributes to $1$
occupied state.  The Hamiltonian operator is represented in a planewave
basis set.   All the examples presented in this section consists of $32$ atoms.
Initially, the atoms are at their equilibrium positions, and
the distance between each atom and its nearest neighbor is set to $10$
a.u..
Starting from the equilibrium position, each ion is given a
finite velocity so that the velocity on the centroid of mass is $0$.
In the numerical experiments below, the system contains only one single phonon,
which is obtained by assigning an initial velocity $v_{0}\propto
(1,-1,1,-1,\cdots)$ to the atoms.  We denote by $\Omega^\text{Ref}$ the
corresponding phonon frequency.  We choose $v_{0}$ so that
$\tfrac{1}{2} mv_0^2 = k_B T_{\ion}$, where $k_B$ is the Boltzmann
constant and $T_{\ion}$ is $10$ K to make
sure that the system is in the linear response regime.
In the atomic unit, the
mass of the electron is $1$, and the mass of each nuclei is set to
$42000$.
By adjusting the parameters $\{\sigma_{I}\}$, the 1D DFT model
model can be tuned to resemble an insulating (with $\sigma_{I}=2.0$)
or a metallic system (with $\sigma_{I}=6.0$) throughout the MD simulation.
Fig.~\ref{fig:eigH} shows the spectrum of the
insulating and the metallic system after running $1000$
BOMD steps with converged SCF iteration.

In the linear response regime, we measure the error of the phonon
frequency calculated from TRBOMD.  This can be done in two ways.
The first is given by Eq.~\eqref{eqn:OmegaPerturb},
namely,
all quantities in the big parentheses in Eq.~\eqref{eqn:OmegaPerturb}
can be directly obtained by using the finite difference method at the
equilibrium position $\vR^{*}$.
The second is to
explore the fact that in the linear response regime, there is linear
relation between the force and the atomic position as in
Eq.~\eqref{eqn:LinearForce}, \ie Hooke's law
\begin{equation}
	f_{I}(t_l) \approx  -m \sum_{J}\mc{D}_{IJ} \wt{R}_{J}(t_l)
	\label{}
\end{equation}
holds approximately at each time step.
Here $\{f_{I}(t_l)\}$ and $\{\wt{R}_{I}(t_l)\}$ are obtained from
the trajectory of the TRBOMD simulation directly.  To numerically
compute $\mc{D}_{IJ}$, we solve the least square problem
\begin{equation}
	\min_{\mc{D}} \sum_{l, I} \, \Bigl\lVert f_{I}(t_l) + m \sum_{J}\mc{D}_{IJ}
	\wt{R}_{J}(t_l) \Bigr\rVert^2
	\label{eqn:LS}
\end{equation}
which yields
\begin{equation}
	\mc{D} = -\frac{1}{m}S^{fR}\left(S^{RR}\right)^{-1},
	\label{}
\end{equation}
where
\begin{equation}
	S^{fR}_{IJ} = \sum_{l} f_{I}(t_l)      \wt{R}_{J}(t_l),\quad
	S^{RR}_{IJ} = \sum_{l} \wt{R}_{I}(t_l) \wt{R}_{J}(t_l).
	\label{}
\end{equation}
The frequencies $\{\wt{\Omega}_{l}\}$ can be obtained by diagonalizing the
matrix $\mc{D}$.  Similarly one can perform the calculation for the
accurate BOMD simulation and obtain the exact value of the frequencies
$\{\Omega_{l}\}$.

In order to compare the performance among BOMD, TRBOMD and CPMD,
we define the following relative errors
\begin{align}
\err_\Omega^\text{Hooke} & = \frac{\wt{\Omega}^\text{Hooke} - \Omega^\text{Ref}}{\Omega^\text{Ref}},  \\
\err_\Omega^\text{LR} & = \frac{\wt{\Omega}^\text{LR} - \Omega^\text{Ref} }{\Omega^\text{Ref}}, \\
\err_{\wb{E}} & = \frac{\wb{E}-\wb{E}^\text{Ref}}{\wb{E}^\text{Ref}}, \\
\err_R^{L^2} & = \frac{\|R_1(t) - R_1^\text{Ref}(t)\|_{L^2}}{\|R_1^\text{Ref}(t)\|_{L^2}}, \\
\err_R^{L^\infty} & = \frac{\|R_1(t) - R_1^\text{Ref}(t)\|_{L^\infty}}{\|R_1^\text{Ref}(t)\|_{L^\infty}},
\end{align}
where the results from BOMD with convergent SCF iteration are taken to be
corresponding reference values, $\wb{E}$ is the average total energy over time,
the frequencies $\wt{\Omega}^\text{Hooke}$ and $\Omega^\text{Ref}$ are obtained via solving
the least square problem \eqref{eqn:LS},
the frequency $\wt{\Omega}^\text{LR}$ is measured by Eq.~\eqref{eqn:OmegaPerturb}
with finite difference methods,
and $R_1(t)$ is the trajectory of the left most atom.

\subsection{Numerical comparison between BOMD and TRBOMD}

\begin{figure}
\centering
\subfigure[Insulator.]
{\includegraphics[width=8cm,height=6cm]{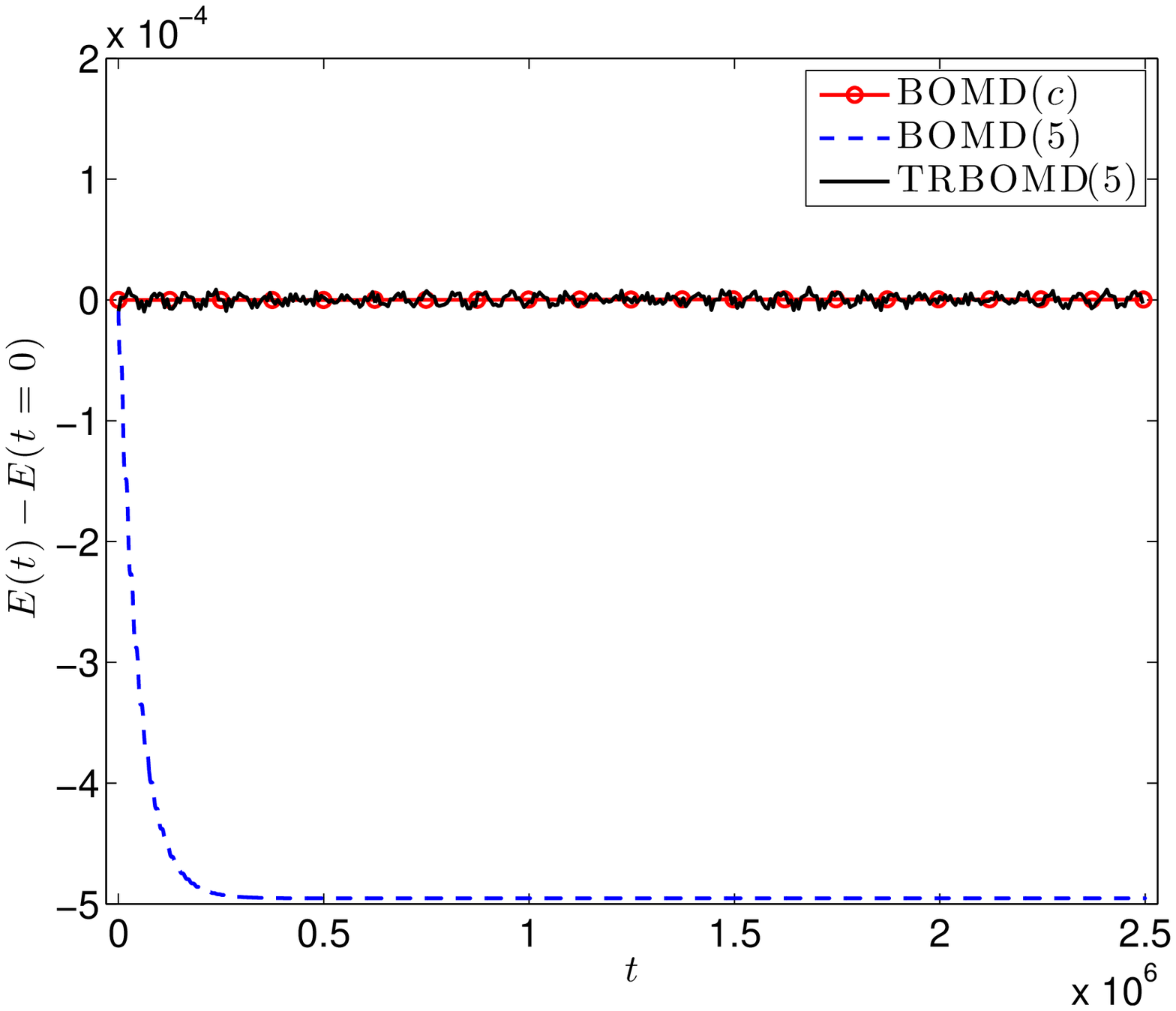}}
\subfigure[Metal.]
{\includegraphics[width=8cm,height=6cm]{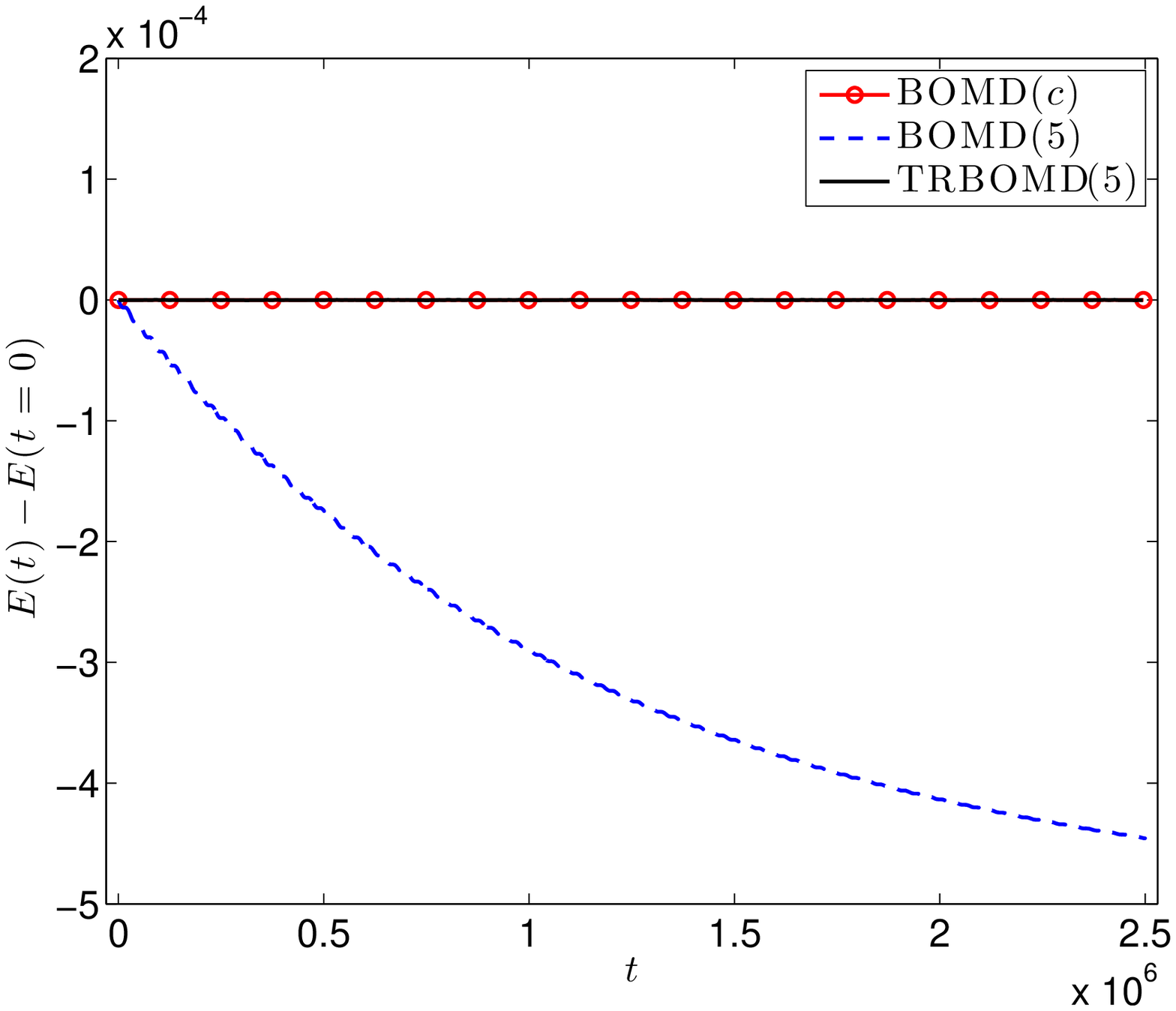}}
\caption{{\small The energy fluctuations around the starting energy $E(t=0)$ as a function of time.
The time step is $\Delta t=250$, the final time is $2.50\E$+$06$ and $\omega = 1/\Delta t=4.00\E$-$03$.
The simple mixing with the Kerker preconditioner is applied in SCF cycles.
BOMD($c$) denotes the BOMD simulation with converged SCF iteration,
and BOMD($n$) (resp. TRBOMD($n$)) represents the BOMD (resp. TRBOMD) simulation with $n$ SCF iterations per time step.
It shows clearly that BOMD($5$) produces
large drift for both insulator (a) and metal (b),
but TRBOMD($5$) does not.
}}\label{drift-compare}
\end{figure}

\begin{figure}
\centering
\subfigure[Insulator.]
{\includegraphics[width=8cm,height=6cm]{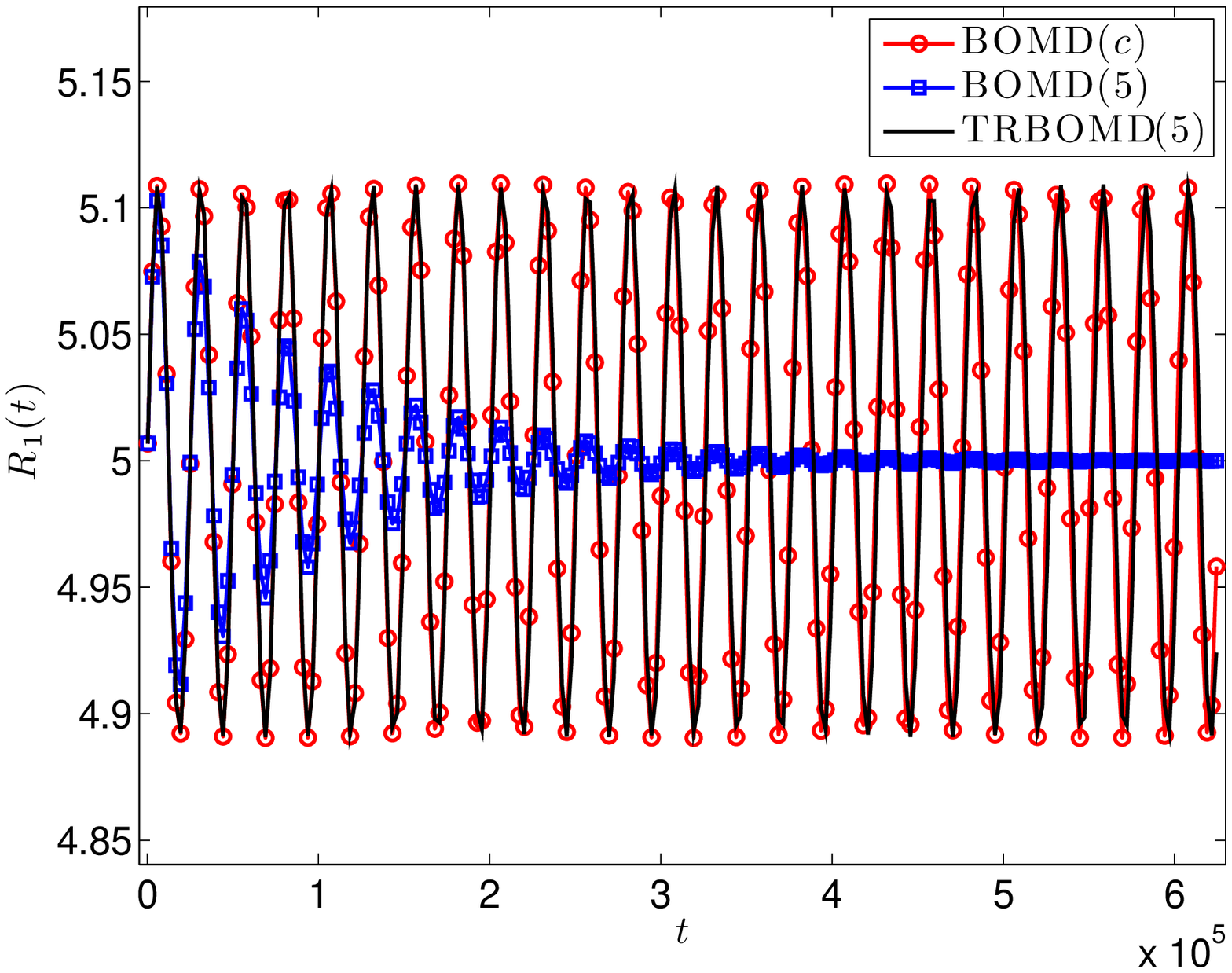}}
\subfigure[Metal.]
{\includegraphics[width=8cm,height=6cm]{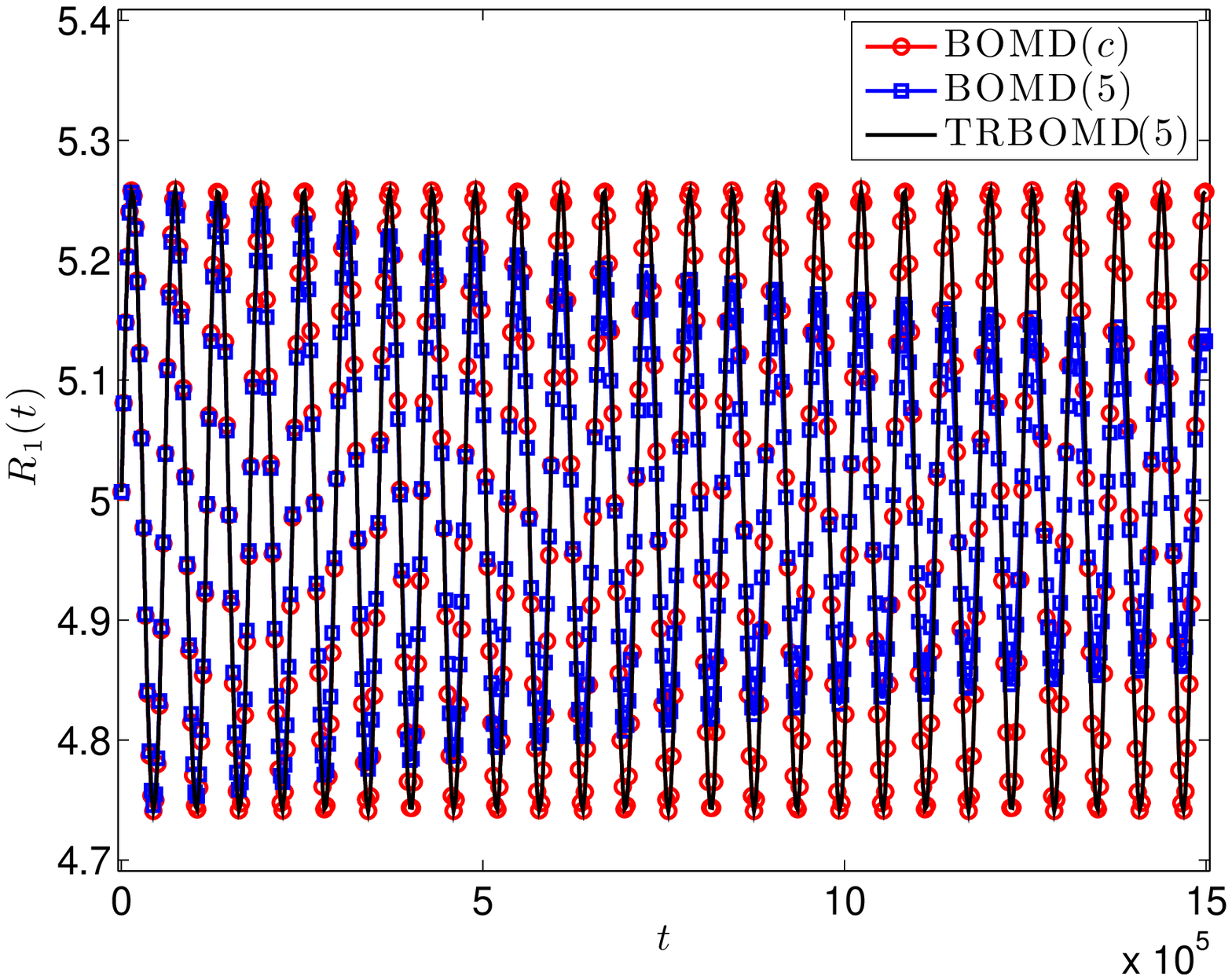}}
\caption{{\small The position of the left most atom as a function of time.
The settings are the same as those in Fig.~\ref{drift-compare}.
It shows clearly that the trajectory from TRBOMD($5$) almost coincides with that from BOMD($c$).
However, for BOMD($5$),  the atom will cease oscillation after a while.
}}\label{position-compare}
\end{figure}

The first run is to validate the performance of TRBOMD.
We set the time step $\Delta t=250$, the artificial frequency $\omega=\frac{1}{\Delta t}=4.00\E$-$03$,
the final time $T=2.50\E$+$06$ and employ the simple mixing with step length $\alpha=0.3$
and the Kerker preconditioner in SCF cycles.
Fig.~\ref{drift-compare} plots the energy drift for BOMD with the converged SCF iteration (denoted by BOMD($c$))
where the tolerance is $1.00\E$-$08$,
BOMD with $5$ SCF iterations per time step (denoted by BOMD($5$)) and
TRBOMD with $5$ SCF iterations per time step (denoted by TRBOMD($5$)).
We see clearly there that BOMD($5$) produces
large drift for both insulator and metal,
but TRBOMD($5$) does not.
Actually, from Table~\ref{tab:error},
the relative error in the average total energy over time between TRBOMD($5$) and BOMD($c$)
is under $1.30\E$-$05$, but
BOMD($c$) needs about average $45$ SCF iterations
per time step to reach the tolerance $1.00\E$-$08$.
Fig.~\ref{position-compare} plots corresponding
trajectory of the left most atom during about the first $25$ periods
and shows that
the trajectory from TRBOMD($5$) almost coincides with that from BOMD($c$),
which is also confirmed by
the data of $\err_R^{L^2}$ and $\err_R^{L^\infty}$ in Table~\ref{tab:error}.
However, for BOMD($5$),
the atom will cease oscillation after a while.
A similar phenomena occurs for other atoms.  In Table~\ref{tab:error},
we present more results for TRBOMD($n$) with $n=3,5,7$. We observe there that
TRBOMD($n$) gives more accurate results with larger $n$, and
$\err_\Omega^\text{Hooke}$ has a similar behavior
as $n$ increases to $\err_\Omega^\text{LR}$, which is in accord with our
previous linear response analysis in Sec.~\ref{sec:lr}.

\begin{table}[h]
  \centering
  \caption{{\small The errors for TRBOMD($n$).
  The settings are the same as those in Fig.~\ref{drift-compare} except for the number of SCF iterations.}}
  \label{tab:error}
\begin{tabular}{ccccccc}
\toprule
\multicolumn{6}{c}{Insulator: $\Omega^\text{Ref}=2.51\E$-$04$, $\wb{E}^\text{Ref}=8.66\E$-$01$}  \\
\midrule
$n$ & $\err_\Omega^\text{LR}$ & $\err_\Omega^\text{Hooke}$ & $\err_{\wb{{E}}}$ & $\err_R^{L^2}$ & $\err_R^{L^\infty}$ \\ \midrule
$3$	& $-6.53\E$-$03$          & $-1.63\E$-$02$          & $-7.63\E$-$05$ &	$2.26\E$-$02$ &	$4.25\E$-$02$  \\ \midrule
$5$	& $-1.08\E$-$03$          &	$-2.38\E$-$03$          & $-1.30\E$-$05$ &	$1.27\E$-$02$ &	$2.92\E$-$02$  \\ \midrule
$7$	& $-2.76\E$-$04$          &	$-5.41\E$-$04$          & $-3.32\E$-$06$	&   $3.02\E$-$03$ &	$7.22\E$-$03$ \\ \midrule\midrule
\multicolumn{6}{c}{Metal: $\Omega^\text{Ref}=1.06\E$-$04$, $\wb{E}^\text{Ref}=5.28\E$-$01$}  \\
\midrule
$3$ & $-2.65\E$-$04$          &	$-6.92\E$-$04$          & $-4.36\E$-$06$ &	$3.86\E$-$03$ &	$8.95\E$-$03$ \\ \midrule
$5$	& $-3.65\E$-$05$          & $-7.31\E$-$05$          & $-4.44\E$-$07$ &	$4.14\E$-$04$ &	$9.60\E$-$04$ \\ \midrule
$7$ & $-5.24\E$-$06$          & $\;\;\;2.93\E$-$06$	          & $-1.10\E$-$07$ &	$1.63\E$-$05$ &	$3.78\E$-$05$ \\
\bottomrule
\end{tabular}
\end{table}

\begin{figure}
\centering
\subfigure[Insulator.]
{\includegraphics[width=8cm,height=6cm]{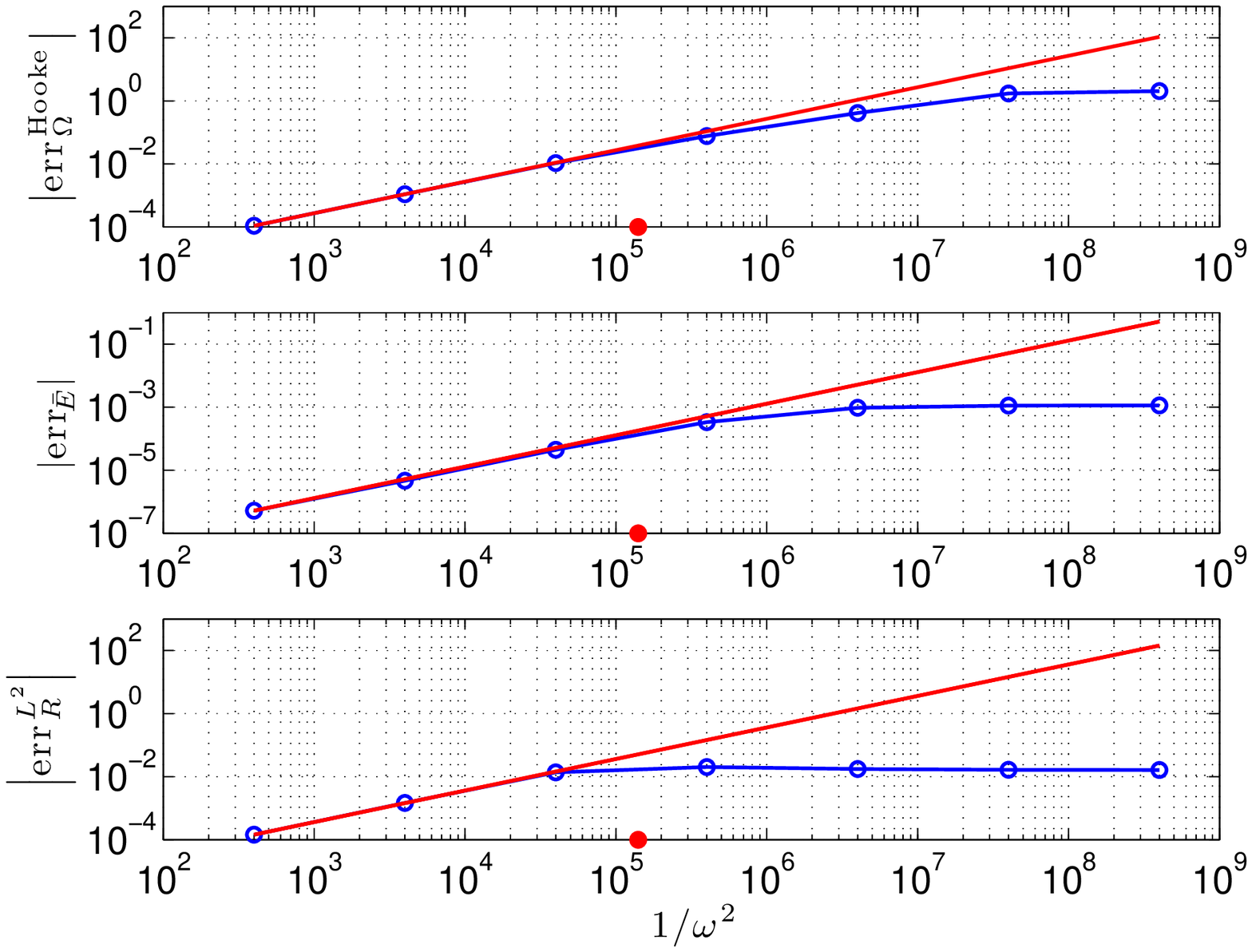}}
\subfigure[Metal.]
{\includegraphics[width=8cm,height=6cm]{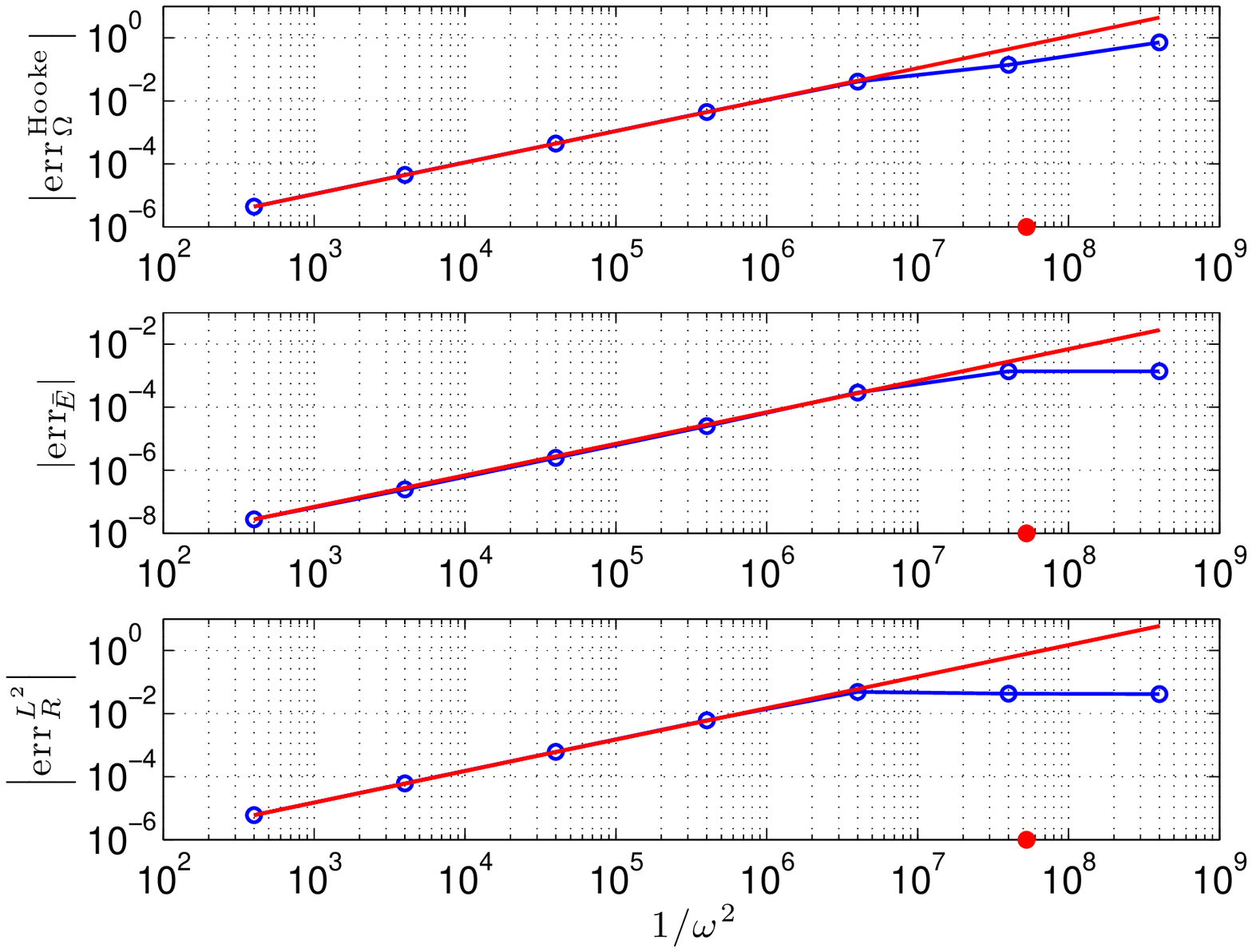}}
\caption{{\small The absolute value of the error
for TRBOMD($3$) as a function of $1/\omega^2$ in logarithmic scales.
The time step is $\Delta t=20$ and the final time is $6.00\E$+$05$.
For the readers' reference, within each plot,
the red straight line denotes corresponding linear dependence
while the red solid point in $x$ axis represents the critical value of ${\lambda_{\min}(\mc{K})}/{\lambda_{\max}(\mc{D})}$.
}}\label{rate-compare}
\end{figure}

According to Eq.~\eqref{eqn:OmegaPerturb}, we have that
$\err_\Omega^\text{LR}$ is proportional to ${1}/{\omega^2}$ for large
$\omega$.  We verify this behavior using TRBOMD($3$) as an example.
In this example, a smaller time step $\Delta t=20$ is set to allow
bigger artificial frequency $\omega$, the final time is
$T=6.00\E$+$05$, and the simple mixing with $\alpha=0.3$ and the
Kerker preconditioner is applied in SCF iterations.  For TRBOMD($3$)
under these settings, we have
$\lambda_{\min}(\mc{K})\simeq8.81\E$-$03$ for the insulator and
$\lambda_{\min}(\mc{K})\simeq5.92\E$-$01$ for the metal, and thus the
critical values of ${(\Omega^\text{Ref})^2}/{\lambda_{\min}(\mc{K})}$
in Eq.~\eqref{eqn:stable1} are about $7.12\E$-$06$ and $1.90\E$-$08$,
respectively.  We choose $\omega^2=2.50\E$-$03$, $2.50\E$-$04$,
$2.50\E$-$05$, $2.50\E$-$06$, $2.50\E$-$07$, $2.50\E$-$08$,
$2.50\E$-$09$, and plot in Fig.~\ref{rate-compare} the absolute values
of $\err_\Omega^\text{Hooke}$, $\err_{\wb{E}}$, $\err_R^{L^2}$ for
TRBOMD($3$) as a function of $1/\omega^2$ in logarithmic scales.  When
$1/\omega^2 \ll {\lambda_{\min}(\mc{K})}/{(\Omega^\text{Ref})^2}$,
Fig.~\ref{rate-compare} shows clearly that all of
$|\err_\Omega^\text{Hooke}|$, $|\err_{\wb{E}}|$, $|\err_R^{L^2}|$
depend linearly on $1/\omega^2$.  The error $\err_R^{L^\infty}$ has a
similar behavior to $\err_R^{L^2}$ and is skipped here for saving
space.

\begin{figure}
\centering
\subfigure[The energy drift.]
{\includegraphics[width=8cm,height=6cm]{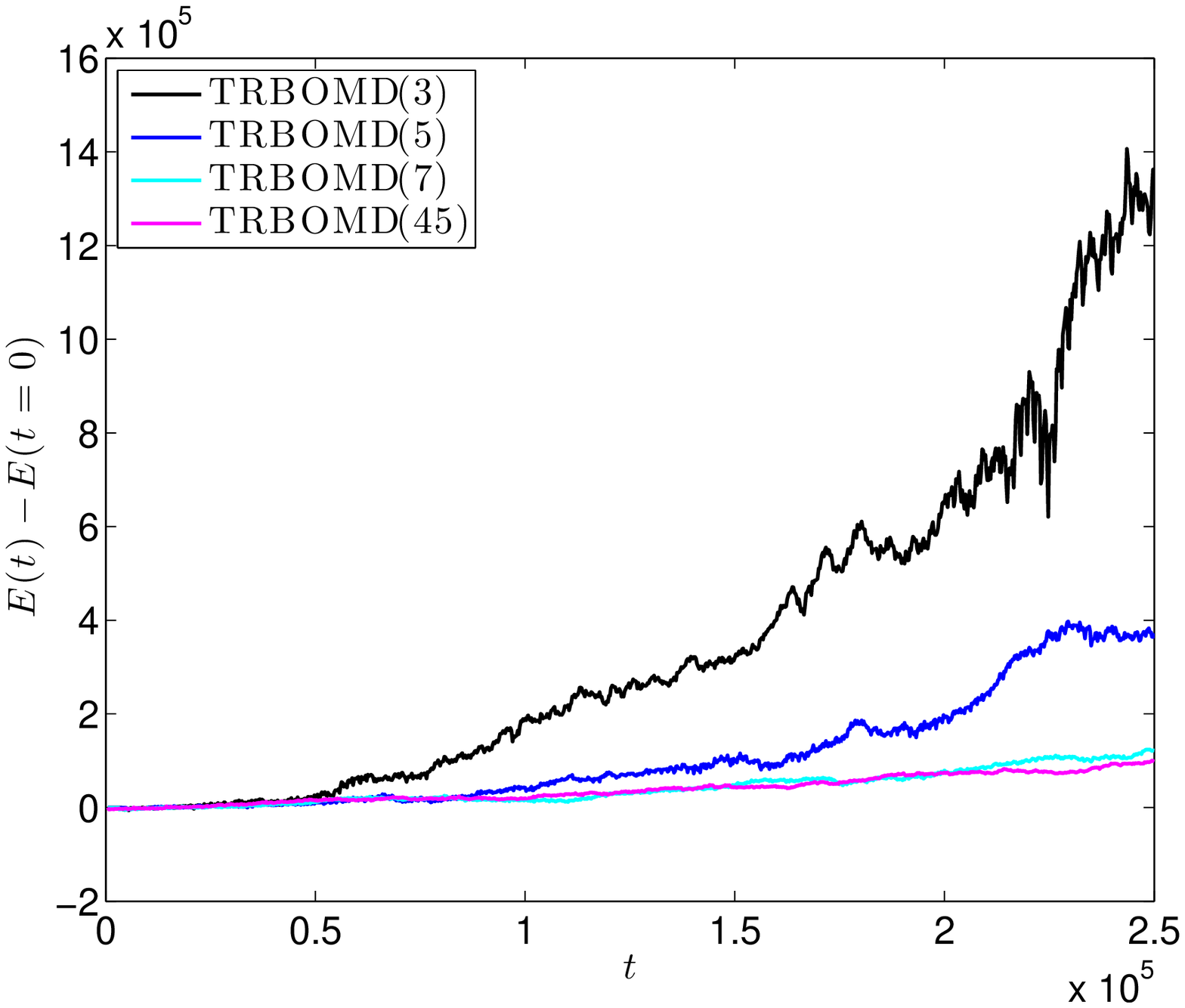}}
\subfigure[The trajectory of the left most atom.]
{\includegraphics[width=8cm,height=6cm]{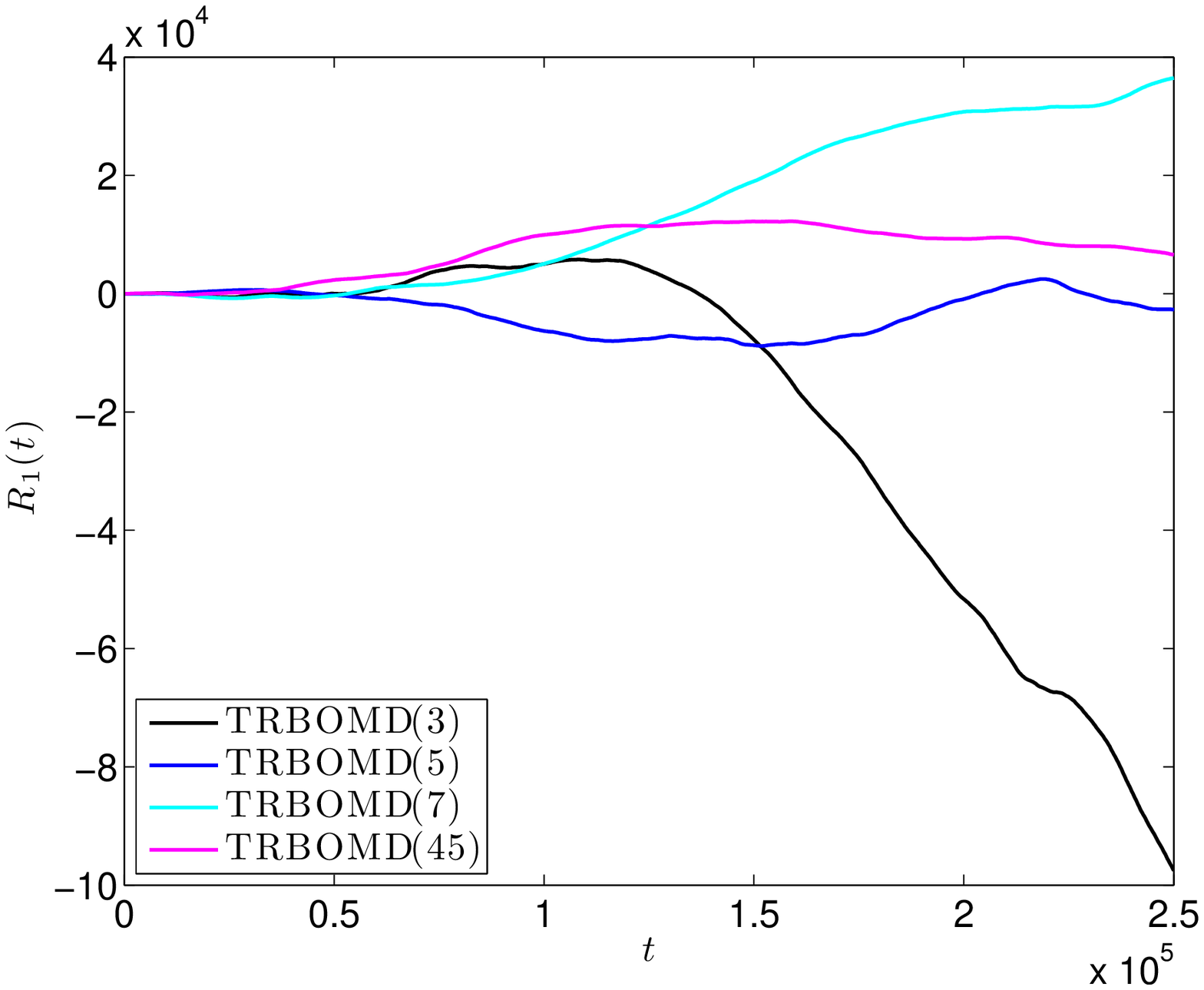}}
\caption{{\small The unstable behavior of TRBOMD with the simple mixing for the insulator.
The time step is $\Delta t=250$, the final time is $2.50\E$+$05$ and $\omega = 1/\Delta t=4.00\E$-$03$.
}}\label{unstable}
\end{figure}

The last example illustrates the possible unstable behavior of TRBOMD
when the stability condition~$\lambda_{\min}(\mc{K}) > 0$ in
Eq.~\eqref{eqn:stability} is violated.
Here we take the insulator as an example and set the time step $\Delta
t=250$, the final time to $2.50\E$+$05$, and the artificial frequency
$\omega=\frac{1}{\Delta t}=4.00\E$-$03$.  The simple mixing with
$\alpha=0.3$ is now applied in SCF iterations.  Under these setting,
we have $\lambda_{\min}(\mc{K})<0$, \eg
$\lambda_{\min}(\mc{K})=-2.42\E$+$03$ for TRBOMD($3$).
Fig.~\ref{unstable}(a) plots the energy drift for TRBOMD($n$) with
$n=3,5,7,45$. We see clearly there that TRBOMD is unstable even using
$45$ SCF iterations per time step (recall that BOMD($c$) in the first
run needs about average $45$ SCF iterations per time
step). Fig.~\ref{unstable}(b) plots corresponding trajectory of the
left most atom and shows that the atom is driven wildly by the
non-convergent SCF iteration.

\subsection{Numerical comparison between TRBOMD and CPMD}

We now present some numerical examples for CPMD illustrating the
difference between CPMD and TRBOMD. As we have discussed, TRBOMD is
applicable to both metallic and insulting systems, while CPMD becomes
inaccurate when the gap vanishes.
To make this statement more concrete, we apply CPMD to the same atom chain system.
We implement CPMD using standard velocity Verlet scheme combined with
RATTLE for the orthonormality constraints~\cite{RyckaertCiccottiBerendsen1977, CiccottiFerrarioRyckaert1982, Andersen1983}.

\begin{figure}[htp]
\centering
\subfigure[Insulator.]
{\includegraphics[width=8cm]{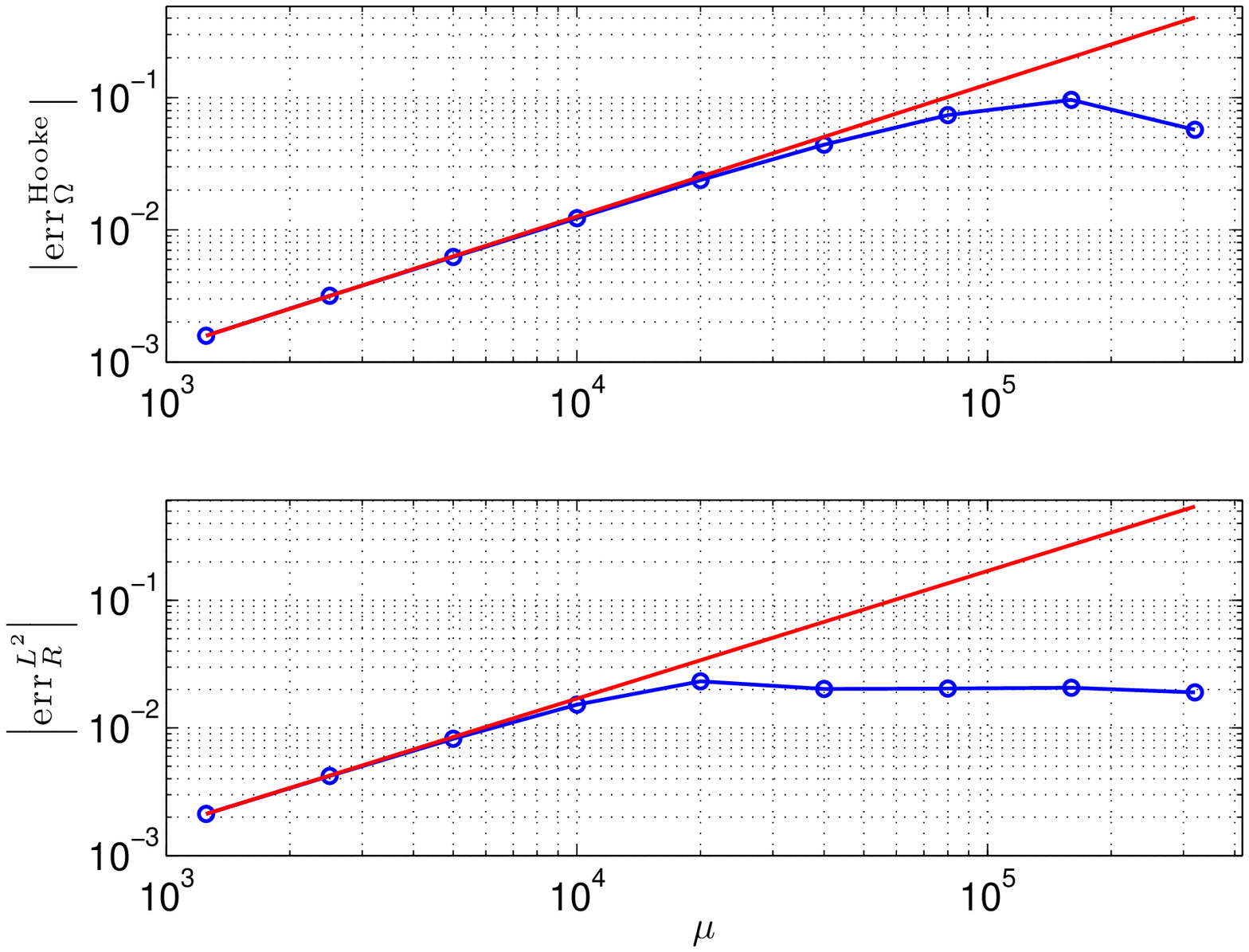}}
\subfigure[Metal.]
{\includegraphics[width=8cm]{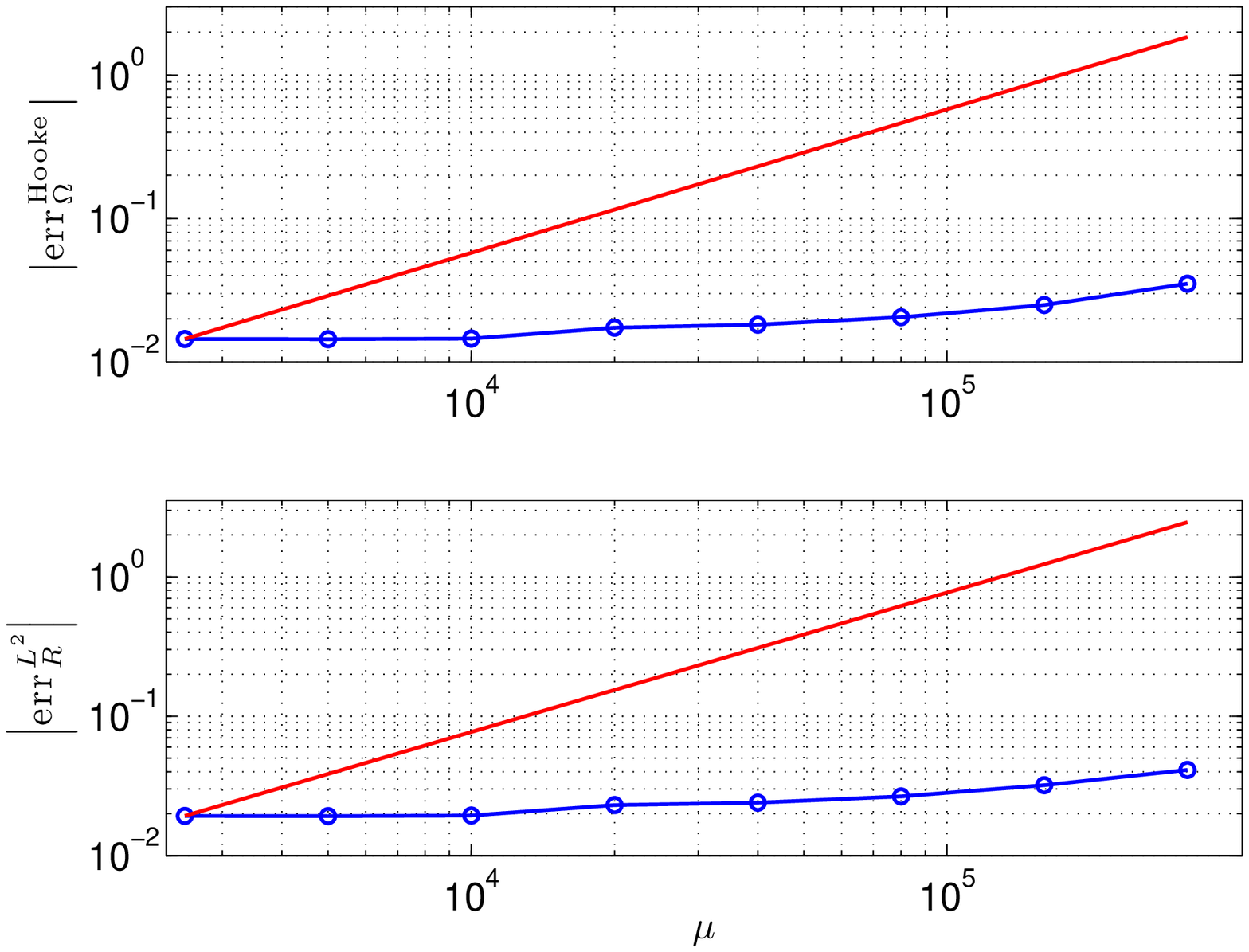}}
\caption{{\small The absolute value of the error
for CPMD as a function of $\mu$ in logarithmic scales.
The time step is $\Delta t=20$ and the final time is $6.00\E$+$05$.}}\label{fig:cpmd-rate}
\end{figure}

We present in Fig.~\ref{fig:cpmd-rate} the error of CPMD simulation for
different choices of fictitious electron mass $\mu$. We study the
relative error of the phonon frequency $\err_\Omega^\text{Hooke}$, the
relative error of position of the left-most atom measured in $L^2$
norm, \ie $\err_R^{L^2}$. We
observe in Fig.~\ref{fig:cpmd-rate}(a) linear convergence of CPMD to
the BOMD result as the parameter $\mu$ decreases. This is consistent
with our analysis. Recall that in CPMD, $\mu$ plays a similar role as
$\omega^{-2}$ in TRBOMD. For the metallic example, the behavior is quite different, actually Fig.~\ref{fig:cpmd-rate}(b) shows a systematic error as $\mu$ decreases. For metallic system, as the spectral gap vanishes, the adiabatic separation between ionic and electronic degrees of freedom cannot be achieved no matter how small $\mu$ is. The adiabatic separation for TRBOMD on the other hand relies on the choice of an effective $\rhoscf$, and hence TRBOMD also works for metallic system as Fig.~\ref{rate-compare} indicates.

\begin{figure}[htp]
\centering
\subfigure[Insulator.]
{\includegraphics[width=8cm]{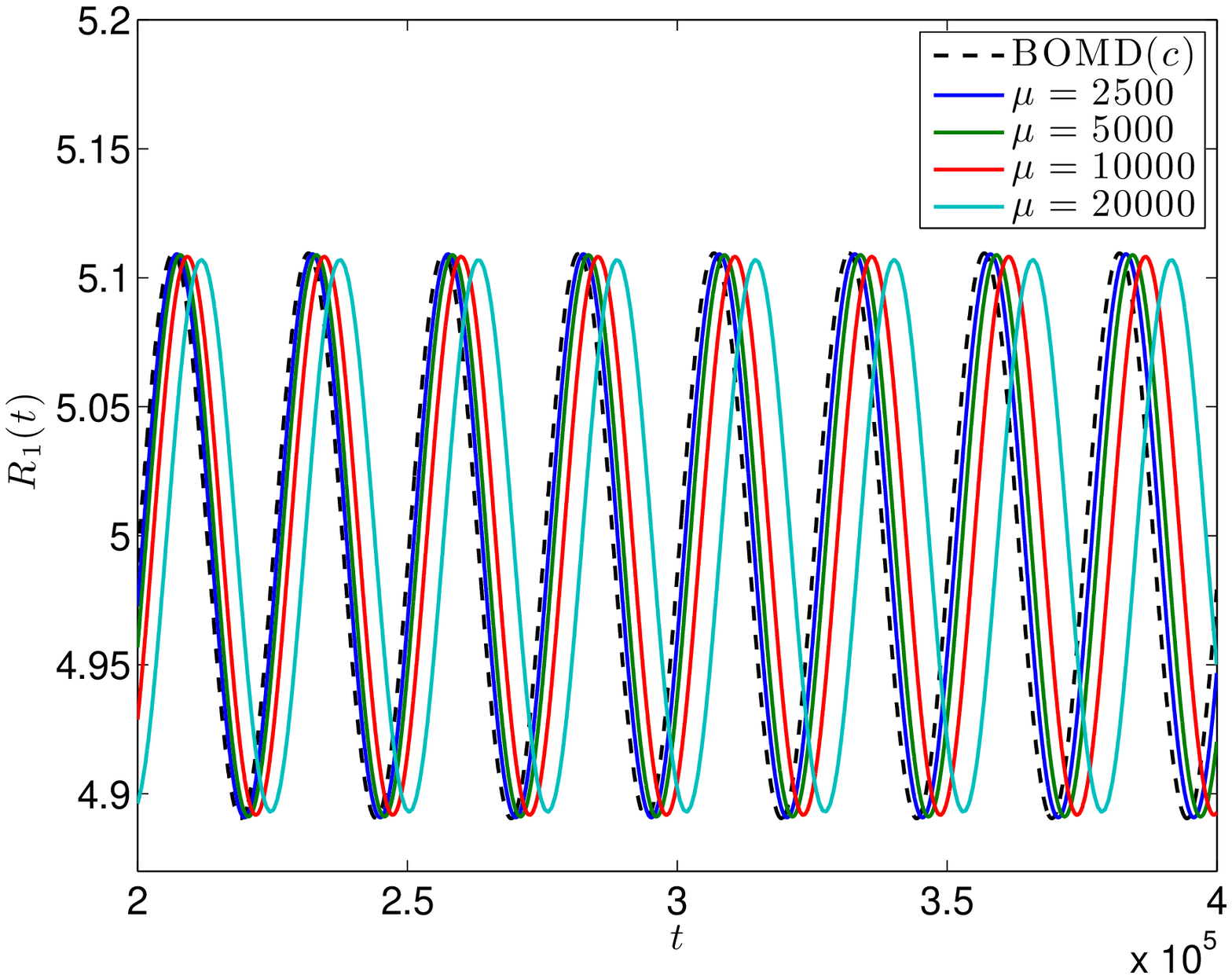}}
\subfigure[Metal.]
{\includegraphics[width=8cm]{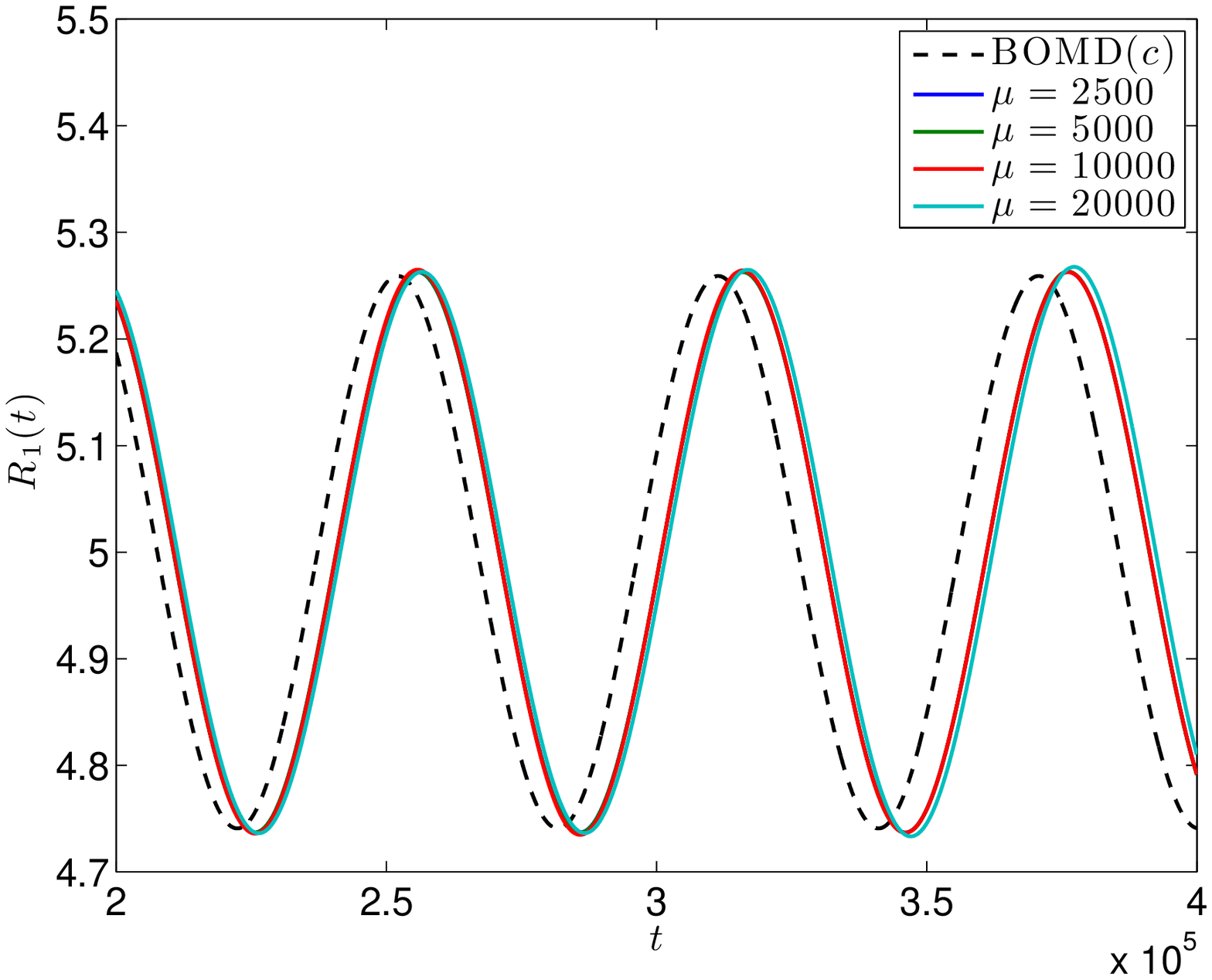}}
\caption{{\small The trajectory of the position of the left-most atom.
Dashed line is the result from BOMD with converged SCF iteration.
Colored solid lines are results from CPMD with fictitious electron mass
$\mu = 2500$, $5000$, $10000$, and $20000$. The time step is $\Delta
t=20$, the trajectory plotted is within the time interval
$[2.00\E$+$05$, $4.00\E$+$05]$}.}
\label{fig:cpmd-traj}
\end{figure}

The different behavior of CPMD for insulating and metallic systems is
further illustrated by Fig.~\ref{fig:cpmd-traj} which shows the
trajectory of the position of the left-most atom during the simulation.
The phase error is apparent from the two subfigures. While the phase
error decreases so that the trajectory approaches that of BOMD for
insulator in Fig.~\ref{fig:cpmd-traj}(a), the result in
Fig.~\ref{fig:cpmd-traj}(b) shows a systematic error for metallic
system.

\section{Beyond the linear response regime: Non-equilibrium dynamics}\label{sec:nonlinear}

The discussion so far has been limited to the linear response regime so
that we can make linear approximations for the degrees of freedom of
both nuclei and electrons. In this case, as the system becomes linear,
explicit error analysis has been given. For practical
applications, we will be also interested in non-equilibrium nuclei
dynamics so that the deviation of atom positions is no longer small. In
this section, we will investigate the non-equilibrium case using
averaging principle (see \eg \cite{Ebook, PavliotisStuart} for general introduction on averaging principle).

\begin{figure}[htb]
  \centering
  \includegraphics[width=8cm]{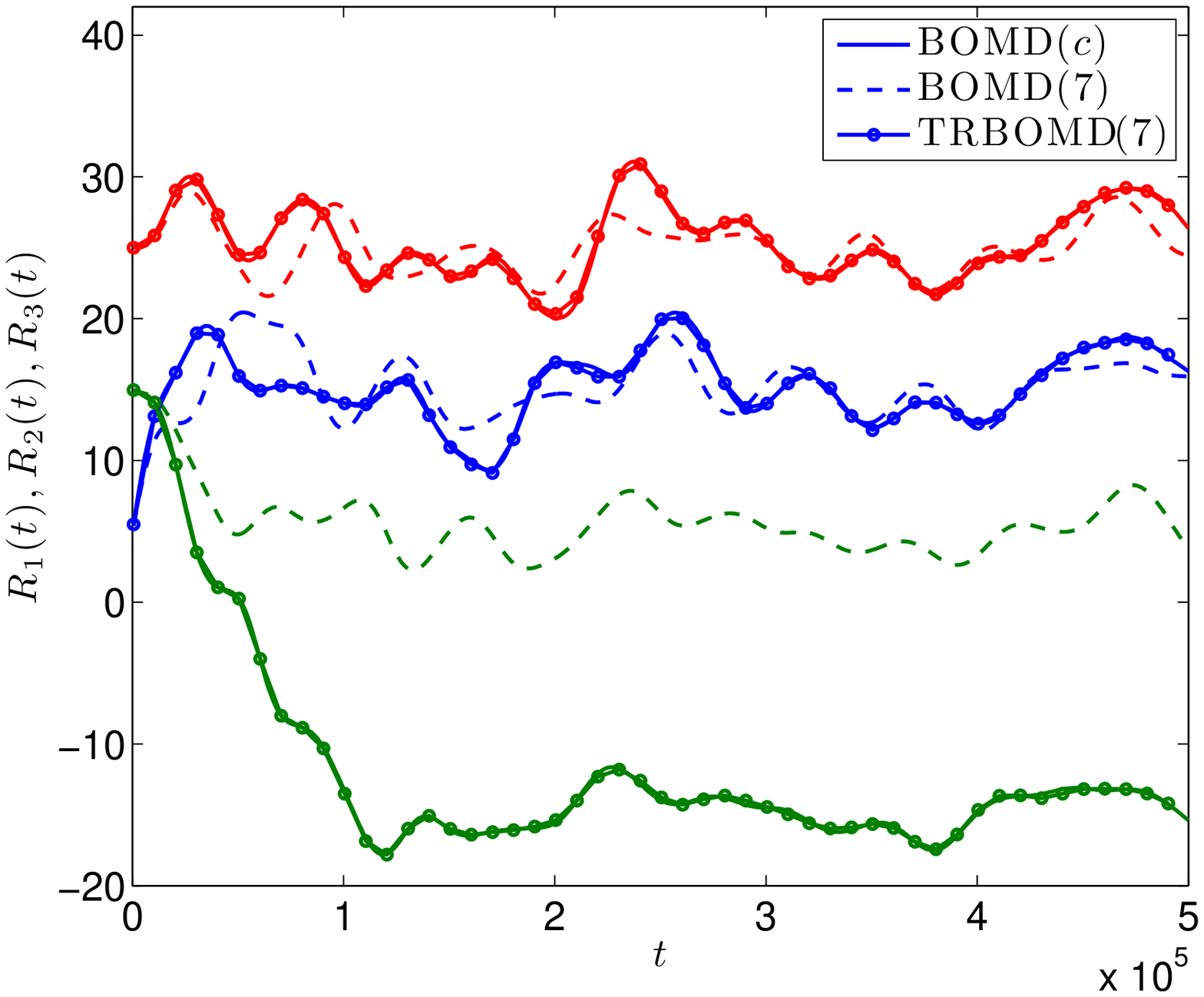}
  \caption{Comparison of trajectories of the first three atoms from the
	left for a non-equilibrium system.
  Different atoms are distinguished by color (blue for the initially left-most atom;
  green for the initially second left-most atom; red for the initially third left-most atom).
  Solid lines are results from BOMD($c$); 
  circled lines are results from TRBOMD($7$); 
  dashed lines are results from BOMD($7$). 
  It is evident that while results from BOMD with a non-convergent SCF iteration have a huge deviation,
  the results from TRBOMD are hardly distinguishable from the ``true'' results from BOMD.
  }\label{fig:R_drift_met}
\end{figure}

\begin{figure}[htbp]
  \centering
  \includegraphics[width=8cm]{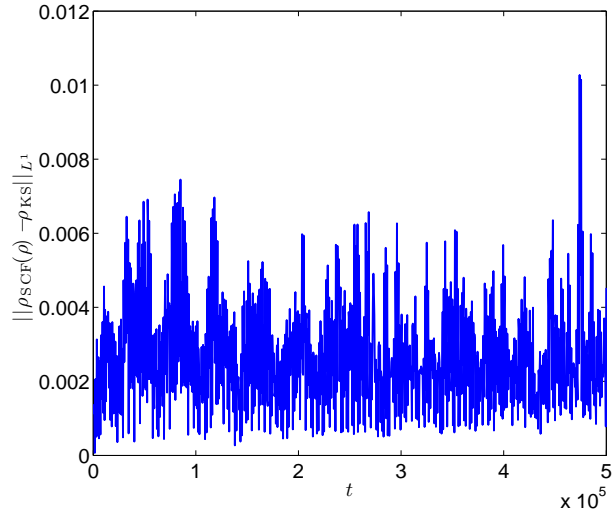}
  \caption{The difference of $\rhoscf$ with the converged electron
	density of SCF iteration (denoted by $\rho_\text{KS}$) measured in $L^1$ norm along the TRBOMD
	simulation for a non-equilibrium system.}\label{fig:rhodist_drift_met}
\end{figure}

Let us first show numerically a non-equilibrium situation for the
atom chain example discussed before. Initially, the $32$ atoms stay at
their equilibrium position. We set the initial velocity so that the
left-most atom has a large velocity towards right and other atoms have
equal velocity towards left. The mean velocity is equal to $0$, so the
center of mass does not move. Fig.~\ref{fig:R_drift_met} shows the
trajectory of positions of the first three atom from the left. We
observe that the results from TRBOMD agree very well with the BOMD
results with convergent SCF iterations. Let us note that in the
simulation, the left-most atom crosses over the second left-most
atom. This happens since in our model, we have taken a $1D$ analog of
Coulomb interaction and the nuclei background charges are smeared out,
and hence the interaction is ``soft'' without hard-core repulsion. In
Fig.~\ref{fig:rhodist_drift_met}, we plot the difference between
$\rhoscf$ and the converged electron density of the SCF iteration (denoted by $\rho_\text{KS}$)
along the TRBOMD simulation. We see that the electron density used in
TRBOMD stays close to the ground state electron density corresponds to
the atom configuration.

To understand the performance of TRBOMD, recall that the equations of motion are given by
\begin{align*}
  m \ddot{R}_I(t) & = - \int \rhoscf(x; \vR(t), \rho(t)) \frac{\partial V_{\ion}(x; \vR(t))}{\partial R_I} \ud x, \\
  \ddot{\rho}(x, t) & = \omega^2 ( \rhoscf(x; \vR(t), \rho(t)) - \rho(x, t))
\end{align*}
To satisfy the adiabatic condition~\eqref{eqn:stable1} from the linear analysis,
$\omega$ here is a large parameter. As a result, the time scales of the motions of the nuclei and of the electrons are quite different: The electronic degrees of freedom move much faster than the nuclear degrees of freedom.

Let us consider the limit $\omega \to \infty$. In this case, we may freeze the $\vR$ degree of freedom in the equation of motion for $\rho$, as $\rho$ changes on a much faster time scale. To capture the two time scale behavior, we introduce a heuristic two-scale asymptotic expansion with faster time variable given by $\tau = \omega t$ (with some abuse of notations):
\begin{equation}
  R(t) = R(t) \quad \text{and} \quad \rho(x, t) = \rho(x, t, \tau),
\end{equation}
and hence
\begin{equation}
  \ddot{\rho}(x, t) = \omega^2 \partial_{\tau}^2 \rho(x, t, \tau) + 2 \omega \partial_{\tau} \partial_t \rho(x, t, \tau) + \partial_t^2 \rho(x, t, \tau).
\end{equation}
Therefore, to the leading order, after neglecting terms of $\Or(\omega^{-1})$, we obtain
\begin{align}
  m \ddot{R}_I(t) & = - \int \rhoscf(x; \vR(t), \rho(t, \tau)) \frac{\partial V_{\ion}(x; \vR(t))}{\partial R_I} \ud x, \label{eqn:twoscaleR}\\
  \partial_{\tau}^2 \rho(x, t, \tau) & = \rhoscf(x; \vR(t), \rho(t, \tau)) - \rho(x, t, \tau). \label{eqn:twoscalerho}
\end{align}
For the equation of motion for $\rho$, note that as $\vR$ only depends on $t$,
the nuclear positions are fixed parameters in Eq.~\eqref{eqn:twoscalerho}.

To proceed, we consider the scenario that $\rho(t, \tau)$ is close to
the ground state electron density corresponding to the current atom
configuration $\rho^{\ast}(\vR(t))$. We have seen from numerical
examples (Fig.~\ref{fig:rhodist_drift_met}) that this is indeed the
case for a good choice of SCF iteration, while we do not have a proof of
this in the general case. Hence, we linearize the map
$\rhoscf$.
\begin{equation}
  \rhoscf(x; \vR, \rho) = \rho^{\ast}(x; \vR) + \int \frac{\delta
	\rhoscf}{\delta \rho}(x,y;\vR, \rho^{\ast}(\vR)) (\rho(y) -
	\rho^{\ast}(y; \vR)) \ud y,
\end{equation}
and Eq.~\eqref{eqn:twoscalerho} becomes
\begin{equation}\label{eq:linearrho}
	\partial_{\tau}^2 \rho(x, t, \tau) = - \mc{K}(\vR) (\rho(x, t, \tau) - \rho^{\ast}(x; \vR(t)))
\end{equation}
where $\mc{K}(\vR)$ is the same as in \eqref{eqn:Kkernel} except it is now defined for each atom configuration $\vR$. Let us emphasize that
here we have only taken the linear approximation for the electronic
degrees of freedom, while keeping the possibly nonlinear dynamics of
$\vR$. This is different from the linear response regime considered
before, where the nuclei motion is also linearized.

Under the stability condition~\eqref{eqn:stability},   
it is easy to
see that for $\rho(t, \tau)$ satisfying Eq.~\eqref{eq:linearrho}, the limit
of time average
\begin{equation}\label{eqn:limitrho}
  \begin{aligned}
    \wb{\rho}(x; \vR(t)) & = \lim_{T \to \infty} \frac{1}{T} \int_0^T
    \rhoscf(x; \vR(t), \rho(t, \tau)) \ud \tau, \\
    & \approx \rho^{\ast}(x; \vR(t)) + \int \frac{\delta \rhoscf}{\delta
      \rho}(x,y;\vR, \rho^{\ast}(\vR)) \biggl( \lim_{T \to \infty}
    \frac{1}{T} \int_0^T \rho(y;t, \tau) - \rho^{\ast}(y;\vR(t)) \ud \tau
    \biggr) \ud y\\
    & = \rho^{\ast}(x; \vR(t)).
  \end{aligned}
\end{equation}
Take the average of Eq.~\eqref{eqn:twoscaleR} in $\tau$, we have
\begin{equation}
  m \ddot{R}_I(t) = - \int \wb{\rho}(x; \vR(t)) \frac{\partial V_{\ion}(x; \vR(t))}{\partial R_I} \ud x.
\end{equation}
Because of Eq.~\eqref{eqn:limitrho}, the above dynamics is given by
\begin{equation}
  m \ddot{R}_I(t) = - \int \rho^{\ast}(x; \vR(t)) \frac{\partial V_{\ion}(x; \vR(t))}{\partial R_I} \ud x
\end{equation}
which agrees with the equation of motion of atoms in BOMD. As we have neglected $\Or(\omega^{-1})$ terms in the averaging, the difference in trajectory of BOMD and TRBOMD is on the order of $\Or(\omega^{-1})$ for finite $\omega$.

\begin{remark}
If we do not make the linear approximation for the electronic degree of freedom, as the map $\rhoscf$ is quite nonlinear and complicated, the analysis of the long time (in $\tau$) behavior of Eq.~\eqref{eqn:twoscalerho} is not as straightforward. In particular, it is not clear to us whether the limit
\begin{equation}
    \wb{\rho}(x; \vR(t)) = \lim_{T \to \infty} \frac{1}{T} \int_0^T
    \rhoscf(x; \vR(t), \rho(t, \tau)) \ud \tau
\end{equation}
exists or how close the limit is to $\rho^{\ast}(x; \vR(t))$ in a fully
nonlinear regime.  One particular difficulty lies in the fact that
unlike BOMD or CPMD, we do not have a conserved Lagrangian for the
TRBOMD. Actually, it is easy to construct much simplified analog of
Eq.~\eqref{eqn:twoscalerho} that the average is different from
$\rho^{\ast}$. For example, if we consider the following analog which only has one degree of freedom $\xi$
\begin{equation}
  \ddot{\xi} = (\xi/2 + a \xi^2) - \xi,
\end{equation}
where $(\xi / 2 + a \xi^2)$ is the analog of $\rhoscf$ here and $a>0$ is a small parameter which characterizes the nonlinearity of the map.
Note that
\begin{equation}
  \ddot{\xi} = - \xi / 2 + a \xi^2 = - \partial_{\xi} (\xi^2 / 4 - a \xi^3 / 3 ).
\end{equation}
The motion of $\xi$ is equivalent to a motion of a particle in an
anharmonic potential. It is clear that if initially $\xi(0) \neq 0$,
the long time average of $\xi$ will not be $0$. Furthermore, if
initially, $\xi(0)$ is too large, the orbit is not closed ($\xi$
escapes the well around $\xi = 0$).  If phenomena similar to this
occur for a general $\rhoscf$, then even in the limit $\omega \to
\infty$, there will be a systematic uncontrolled bias between BOMD and
TRBOMD.  This is in contrast with Car-Parrinello molecular dynamics,
which agrees with BOMD in the limit fictitious mass goes to zero ($\mu
\to 0$) if the adiabatic condition holds.

As a result of this discussion, in practice, when we apply TRBOMD to a
particular system, we need to be cautious whether the electronic degree
of freedom remains around the converged Kohn-Sham electron density,
which is not necessarily guaranteed (in contrast to CPMD for systems with gaps).
\end{remark}

\section{Conclusion}\label{sec:conclusion}

The recently developed time reversible Born-Oppenheimer molecular
dynamics (TRBOMD) scheme provides a promising way for reducing the
number of self-consistent field (SCF) iterations in molecular dynamics
simulation.  By introducing auxiliary dynamics to the initial guess of
the SCF iteration, TRBOMD preserves the time-reversibility of the NVE
dynamics both at the continuous and at the discrete level, and
exhibits improved long time stability over the Born-Oppenheimer
molecular dynamics with the same accuracy.  In this paper we analyze
for the first time the accuracy and the stability of the TRBOMD
scheme, and our analysis is verified through numerical experiments
using a one dimensional density functional theory (DFT) model without
exchange correlation potential.  The validity of the stability
condition in TRBOMD is directly associated with the quality
of the SCF iteration procedure. In particular, we demonstrate in the
case when the SCF iteration procedure is not very accurate, the
stability condition can be violated and TRBOMD becomes
unstable.  We also compare TRBOMD with the Car-Parrinello
molecular dynamics (CPMD) scheme. CPMD relies on the adiabatic
evolution of the occupied electron states and therefore CPMD works
better for insulators than for metals.  However, TRBOMD may be
effective for both insulating and metallic systems.  The present study
is restricted to NVE system and to simplified DFT models.  The
performance of TRBOMD for NVT system and for realistic DFT systems
will be our future work.

\section*{Acknowledgments}

This work was partially supported by the Laboratory Directed Research
and Development Program of Lawrence Berkeley National Laboratory under
the U.S. Department of Energy contract number DE-AC02-05CH11231 and
the Scientific Discovery through Advanced Computing
(SciDAC) program funded by U.S. Department of Energy, Office of
Science, Advanced Scientific Computing Research and Basic Energy
Sciences (L. L.), the Alfred P. Sloan Foundation (J. L.), and
the National Natural Science Foundation of China under the grant number
11101011 and the Specialized Research Fund for the Doctoral Program of
Higher Education under the grant number 20110001120112 (S.~S.).

\appendix

\section*{Appendix}\label{sec:appB}

Here we derive the perturbation analysis result in
Eq.~\eqref{eqn:OmegaPerturb}.  When deriving the perturbation analysis
below, we use linear algebra notation and do not distinguish matrices
from operators. We use the linear algebra notation, replace all the
integrals by matrix-vector multiplication, and drop all the dependencies
of the electron degrees of freedom $x$ and $y$.  For instance,
$\mc{K}\wt{\rho}$ should be
understood as $\int \mc{K}(x,y)\wt{\rho}(y) \ud y$. We also denote
$\frac{\partial \rho^{*}}{\partial \vR}(x;\vR^*)$ simply by
$\frac{\partial \rho^{*}}{\partial \vR}$, then
Eq.~\eqref{eqn:linearTRBOMD2} can be rewritten as
\begin{equation}
	\begin{pmatrix}
		\ddot{\wt{\vR}}\\
		\ddot{\wt{\rho}}
	\end{pmatrix}
	= A \begin{pmatrix}
		\wt{\vR}\\
		\wt{\rho}
	\end{pmatrix}
	= \left( A_{0} + \frac{1}{\epsilon} A_{1} \right)
	\begin{pmatrix}
		\wt{\vR}\\
		\wt{\rho}
	\end{pmatrix}.
	\label{}
\end{equation}
Here
\begin{equation}
	A_{1} = \begin{pmatrix}
		0 & 0 \\
		0 & -\mc{K}
	\end{pmatrix}
	\label{eqn:A1}
\end{equation}
is a block diagonal matrix, and
\begin{equation}
	A_{0} = \begin{pmatrix}
		-\mc{D} & \vL \\
		\left(\frac{\partial \rho^{*}}{\partial \vR}\right)^T \mc{D} &
		-\left(\frac{\partial \rho^{*}}{\partial \vR}\right)^T \vL
	\end{pmatrix}
	= \begin{pmatrix}
		\mc{I}\\
		-\left(\frac{\partial \rho^{*}}{\partial \vR}\right)^T
	\end{pmatrix}
	\begin{pmatrix}
		-\mc{D} & \vL
	\end{pmatrix}
	\label{eqn:A0}
\end{equation}
is a rank-$M$ matrix. $\mc{I}$ is a $M\times M$ identity matrix. Now
assume the eigenvalues and eigenvectors of $A$ follows the expansion
\begin{equation}
	\lambda = \lambda_{0} + \epsilon \lambda_{1} + \cdots,\quad
	v = v_{0} + \epsilon v_{1} + \cdots.
	\label{}
\end{equation}
Match the equation up to $\Or(\epsilon)$, and
\begin{subequations}
\begin{align}
	A_{1} v_{0} &= 0,\label{eqn:orderm1}\\
	A_{0} v_{0} + A_{1} v_{1} &= \lambda_{0} v_{0},\label{eqn:order0}\\
	A_{0} v_{1} + A_{1} v_{2} &= \lambda_{0} v_{1} + \lambda_{1}
	v_{0}.\label{eqn:order1}
\end{align}
\end{subequations}
Eq.~\eqref{eqn:orderm1} implies that $v_{0}\in \Ker A_{1}$.
Apply the projection
operator $P_{\Ker A_1}$ to both sides of Eq.~\eqref{eqn:order0}, and use
that $v_{0}=P_{\Ker A_1} v_0$, we have
\begin{equation}
	P_{\Ker A_1} A_{0} P_{\Ker A_1}v_{0} = \lambda_{0} P_{\Ker A_1}v_{0}.
	\label{}
\end{equation}
or
\begin{equation}
  \begin{pmatrix}
		-\mc{D} & 0\\
		0 & 0
	\end{pmatrix} v_{0} = \lambda_{0} v_{0}.
	\label{}
\end{equation}
From the eigen-decomposition of $\mc{D}$ in Eq.~\eqref{eqn:phonon1} we
have $\lambda_{0}=-\Omega^2_{l}$ for some $l=1,\ldots,M$. For a fixed
$l$, the corresponding eigenvector to the $0$-th order is
\begin{equation}
	v_{0}=(\vv_l, \vzero)^T.
	\label{}
\end{equation}
From Eq.~\eqref{eqn:order0} we also have
\begin{equation}
	A_{1} v_{1} = \lambda_{0} v_{0} - A_{0}v_{0}
	= \begin{pmatrix}
		\vzero\\
		-\Omega_{l}^2\left(\frac{\partial \rho^{*}}{\partial \vR}\right)^T
		\vv_{l}
	\end{pmatrix},
	\label{}
\end{equation}
and therefore
\begin{equation}
	v_{1} = \Omega_{l}^2\left(\vzero, \mc{K}^{-1}\left[\left(\frac{\partial
	\rho^{*}}{\partial \vR}\right)^T \vv_l \right]\right)^T
	\label{}
\end{equation}
Finally we apply $v_{0}$ to both sides of
Eq.~\eqref{eqn:order1} we have
\begin{equation}
	\lambda_{1} = (v_{0}, A_{0} v_{1}) - (v_{0},\lambda_{0} v_{1}) =
	\Omega_{l}^2 \vv_{l}^T \vL\left[\mc{K}^{-1}\left[\left(\frac{\partial
	\rho^{*}}{\partial \vR}\right)^T \vv_l \right]  \right].
	\label{}
\end{equation}
Therefore
\begin{equation}
	\lambda = -\Omega_{l}^2 + \epsilon \Omega_{l}^2 \vv_{l}^T \vL\left[\mc{K}^{-1}\left[\left(\frac{\partial
	\rho^{*}}{\partial \vR}\right)^T \vv_l \right]  \right] +
	\Or(\epsilon^2)
	\label{}
\end{equation}
In other words, the phonon frequency $\wt{\Omega}_{l}=\sqrt{-\lambda}$
up to the leading order is
\begin{equation}
	\wt{\Omega}_{l} = \Omega_{l}\left( 1 - \frac{1}{2\omega^2}
	\vv_{l}^T \vL\left[ \mc{K}^{-1}\left[ \left(\frac{\partial
	\rho^{*}}{\partial \vR}\right)^T \vv_{l} \right] \right]
	\right) + \Or(1/\omega^4).
	\label{}
\end{equation}
which is Eq.~\eqref{eqn:OmegaPerturb}.

\bibliographystyle{mdpi}
\makeatletter
\renewcommand\@biblabel[1]{#1. }
\makeatother

\end{document}